\documentclass[11pt]{article}

\usepackage{amsfonts,amssymb}
\usepackage[all, knot]{xy}

\newcommand{\lt}{\triangleleft}

\newcommand{\C}{\mathbb{C}}
\newcommand{\Z}{\mathbb{Z}}

\newcommand{\F}{\mathbb{F}}
\newcommand{\Sym}{\mathbb{S}}
\newcommand{\Hom}{\mbox{Hom}}
\newcommand{\Ad}{\mbox{\rm Ad}}

\newcommand{\ad}{ {\mbox{\rm ad}} }

\newcommand{\GG}{{\mathcal G}}

\newtheorem{theorem}{Theorem}[section]

\newtheorem{lemma}[theorem]{Lemma}
\newtheorem{proposition}[theorem]{Proposition}

\newtheorem{definition}[theorem]{Definition}

\newtheorem{example}[theorem]{Example}

\newtheorem{remark}[theorem]{Remark}

\input{epsf.sty}

\begin{document}

\title{Cohomology of the Adjoint of Hopf Algebras }

\author{
J. Scott Carter\footnote{Supported in part by
NSF Grant DMS \#0301095, \#0603926.}
\\ University of South Alabama
\and
Alissa S. Crans \\
Loyola Marymount University
\and
Mohamed
Elhamdadi
\\ University of South Florida
\and
Masahico Saito\footnote{Supported in part by NSF Grant DMS \#0301089,
 \#0603876.}
\\ University of South Florida
}

\maketitle

\begin{abstract}
A cohomology theory of the adjoint of Hopf 
algebras, via deformations,  is presented by means of diagrammatic
techniques.  Explicit calculations are provided in the cases of
group algebras, function algebras on groups, and the bosonization
of the super line. As applications, solutions to the YBE are given
and quandle cocycles are constructed from groupoid cocycles.

\end{abstract}

\section{Introduction}
Algebraic deformation theory \cite{Gerst} can be used to define
$2$-dimensional cohomology in a wide variety of contexts. This
theory has also been understood diagrammatically
\cite{CCES1,MrSt,MrVo} via PROPs, for example.
In this paper, we use
diagrammatic techniques to define a cohomological deformation of
the adjoint map ${\rm ad}(x \otimes y)= \sum S(y_{(1)}) x y_{(2)}$
in an arbitrary Hopf algebra.
We have
concentrated on the diagrammatic
versions here because
diagrammatics have led to topological invariants
\cite{CJKLS,Kuperberg,ReshTur}, diagrammatic methodology is
prevalent in understanding particle interactions and scattering in
the physics literature, and most importantly kinesthetic intuition
can be used to prove algebraic identities.

The starting point for this calculation is a pair of identities
that the adjoint map satisfies and that are sufficient to
construct Woronowicz's solution \cite{Woro} $R= (1 \otimes {\rm
ad})(\tau \otimes 1)(1 \otimes \Delta)$ to the Yang-Baxter
equation (YBE): $(R \otimes 1)(1 \otimes R)(R \otimes 1)=(1
\otimes R)(R \otimes 1)(1 \otimes R)$. We use deformation theory
to define an extension
$2$-cocycle. Then we show that the
resulting $2$-coboundary map, when composed with the Hochschild
$1$-coboundary map
 is trivial.
A $3$-coboundary is defined via the ``movie move" technology.
Applications of this cohomology theory
include
constructing new solutions
to the YBE by deformations
and constructing quandle cocycles from
groupoid cocycles that arise from this theory.

The paper is organized as follows.
Section~\ref{Pre} reviews the definition of
 Hopf 
 algebras, defines
the adjoint map, and illustrates
Woronowicz's solution to the YBE.
Section~\ref{deformsec} contains the deformation theory.
Section~\ref{Diff} defines the chain groups and differentials in general.
Example calculations in the case of a group algebra, the
function algebra on a group,
 and a calculation of the $1$- and
$2$-dimensional cohomology of the bosonization of the superline
are 
presented in Section~\ref{Examples}. Interestingly, the group
algebra and the function algebra on a group are cohomologically
different. Moreover, the conditions that result when a function on
the group algebra satisfies the cocycle condition coincide with
the definition of groupoid cohomology. This relationship is given
in Section~\ref{gpoidsec}, along with a construction of quandle
$3$-cocycles from groupoid $3$-cocycles.
In Section~\ref{Rmat},
we use the deformation cocycles to construct solutions to the
Yang-Baxter equation.

\subsection{Acknowledgements}  JSC and MS gratefully acknowledge the
 support of the NSF
without which substantial portions of the work would not have been
possible. JSC, ME, and MS have benefited from several detailed
presentations on deformation theory that have been given by J\"{o}rg
Feldvoss. AC acknowledges useful and on-going conversations with
John Baez.

\section{Preliminaries}
\label{Pre}

\begin{figure}[htb]
\begin{center}
\mbox{
\epsfxsize=3.5in
\epsfbox{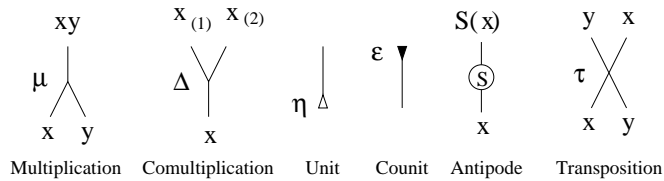}
}
\end{center}
\caption{Operations in Hopf algebras }
\label{Hopfdiag}
\end{figure}

\begin{figure}[htb]
\begin{center}
\mbox{
\epsfxsize=4.5in
\epsfbox{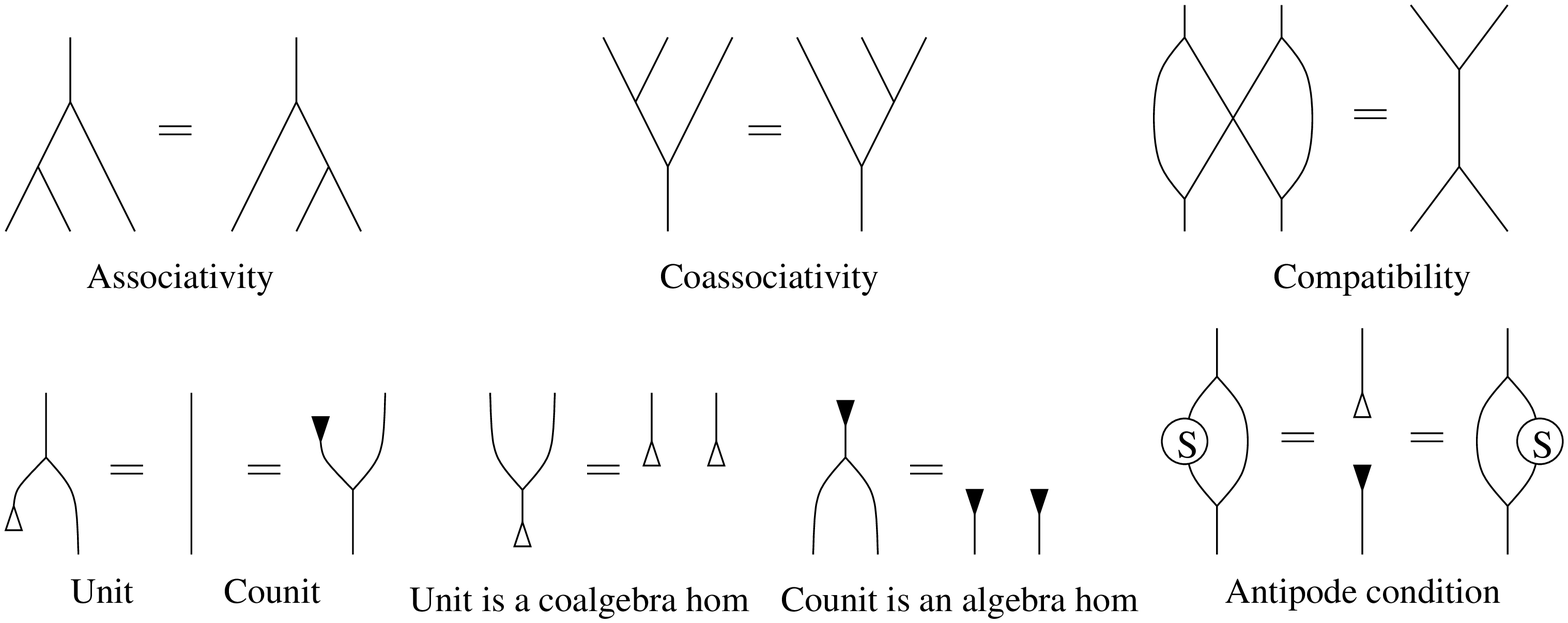}
}
\end{center}
\caption{Axioms of Hopf algebras }
\label{Hopfaxioms}
\end{figure}

We begin by recalling the operations and axioms in Hopf algebras,
and their diagrammatic conventions
depicted in Figures ~\ref{Hopfdiag} and \ref{Hopfaxioms}.

A {\it coalgebra} is a vector space $C$ over a field $k$ together
with a {\it comultiplication} $\Delta: C \rightarrow C \otimes C$
that is bilinear and {\it coassociative}: $(\Delta \otimes 1)
\Delta = (1 \otimes \Delta)\Delta$.  A coalgebra is {\it
cocommutative} if the comultiplication satisfies $\tau \Delta =
\Delta$, where $\tau: C \otimes C \rightarrow  C \otimes C$ is the
transposition $\tau(x \otimes y)=y \otimes x$. A {\it coalgebra
with counit} is a coalgebra with a linear map called the {\it
counit} $\epsilon: C \rightarrow k$ such that $ (\epsilon \otimes
1)\Delta = 1 = (1 \otimes \epsilon)\Delta $ via $k \otimes C \cong
C$.
A  {\it bialgebra} is an algebra $A$ over a field $k$ together
with a linear map called the {\it unit} $\eta: k \rightarrow A$,
satisfying $\eta(a)=a {\bf 1}$ where ${\bf 1} \in A$ is the
multiplicative identity and with an associative multiplication
$\mu: A \otimes A \rightarrow A$ that is also a coalgebra
such
that the comultiplication $\Delta$ is an algebra homomorphism. A
{\it Hopf algebra} is a bialgebra $C$ together with a map called
the {\it antipode} $S: C \rightarrow C$ such that $\mu (S \otimes
1) \Delta = \eta \epsilon = \mu (1 \otimes S) \Delta$, where
$\epsilon$ is the counit.

In diagrams, the compositions of maps are depicted from bottom to
top. Thus a multiplication $\mu$ is represented by a trivalent
vertex with two
bottom
edges representing $A \otimes A$ and one
top
edge representing $A$. Other maps in the definition are
depicted in Fig.~\ref{Hopfdiag} and axioms are depicted in
Fig.~\ref{Hopfaxioms}.

\begin{figure}[htb]
\begin{center}
\mbox{
\epsfxsize=1in
\epsfbox{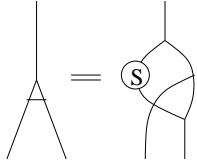}
}
\end{center}
\caption{Adjoint map in a Hopf algebra }
\label{HopfQ}
\end{figure}

Let $H$ be a Hopf algebra. The adjoint map $\Ad_y:  H \rightarrow
H$ for any $y \in H$ is defined by $\Ad_y(x)= S(y_{(1)}) x
y_{(2)}$, where we use the common notation
$\Delta(x)=x_{(1)} \otimes x_{(2)}$ and
$\mu(x \otimes y) = xy$. Its diagram is depicted in
Fig.~\ref{HopfQ}.
Notice the analogy with
group conjugation:
in a group ring
$H=kG$ over a field $k$,  where
$\Delta(y) = y \otimes y$ and $S(y)=y^{-1}$,
we have
$\Ad_y(x )=y^{-1}xy$.

When we view the adjoint map as
a map from $H \otimes H $ to $H$, we use the notation
$$\ad : H \otimes H \rightarrow H , \quad \ad(x \otimes y)=\Ad_y (x) . $$

\begin{figure}[htb]
\begin{center}
\mbox{
\epsfxsize=2.7in
\epsfbox{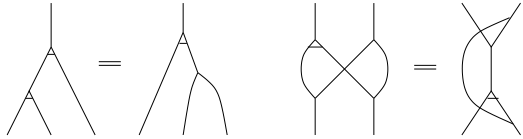}
}
\end{center}
\caption{Conditions for the YBE for Hopf algebras }
\label{Hopfqcond}
\end{figure}

\begin{definition} {\rm
Let $H$ be a Hopf algebra and
$\ad $
 be the adjoint map.
Then the linear  map $R_\ad: H \otimes H \rightarrow H \otimes H $
defined by
$$R_\ad= (1 \otimes \ad)(\tau \otimes 1 )(1 \otimes \Delta)$$
is said to be the $R$-matrix
{\it induced from $\ad$}.
 } \end{definition}

\begin{figure}[htb]
\begin{center}
\mbox{
\epsfxsize=6in
\epsfbox{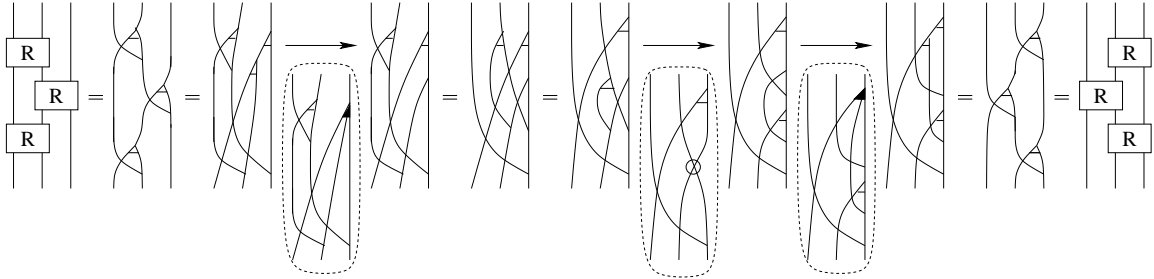}
}
\end{center}
\caption{YBE by the adjoint map }
\label{hopfYBE}
\end{figure}

\begin{lemma}\label{adYBElem}
The $R$-matrix induced from $\ad$ satisfies the YBE.
\end{lemma}
{\it Proof.\/}
In Fig.~\ref{hopfYBE}, it is indicated that
the YBE follows from two properties of the adjoint map:
\begin{eqnarray}
  \ad(\ad  \otimes 1)  &=&  \ad( 1 \otimes \mu)
   \quad {\rm and}  \label{ad1eq} \\
 (\ad \otimes \mu)(1 \otimes \tau \otimes 1)(\Delta
\otimes \Delta) &=& (1 \otimes \mu)(\tau \otimes 1)(1 \otimes
\Delta)(1 \otimes \ad)(\tau \otimes 1)(1 \otimes \Delta).
   \label{ad2eq}
\end{eqnarray}

\begin{figure}[htb]
\begin{center}
\mbox{
\epsfxsize=3.5in
\epsfbox{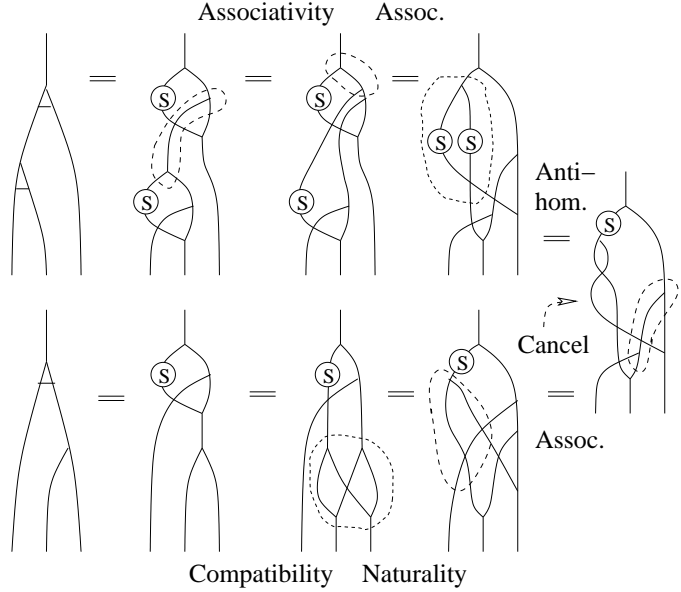}
}
\end{center}
\caption{Adjoint condition $(1)$}
\label{HopfQlemma}
\end{figure}

\begin{figure}[htb]
\begin{center}
\mbox{
\epsfxsize=5.5in
\epsfbox{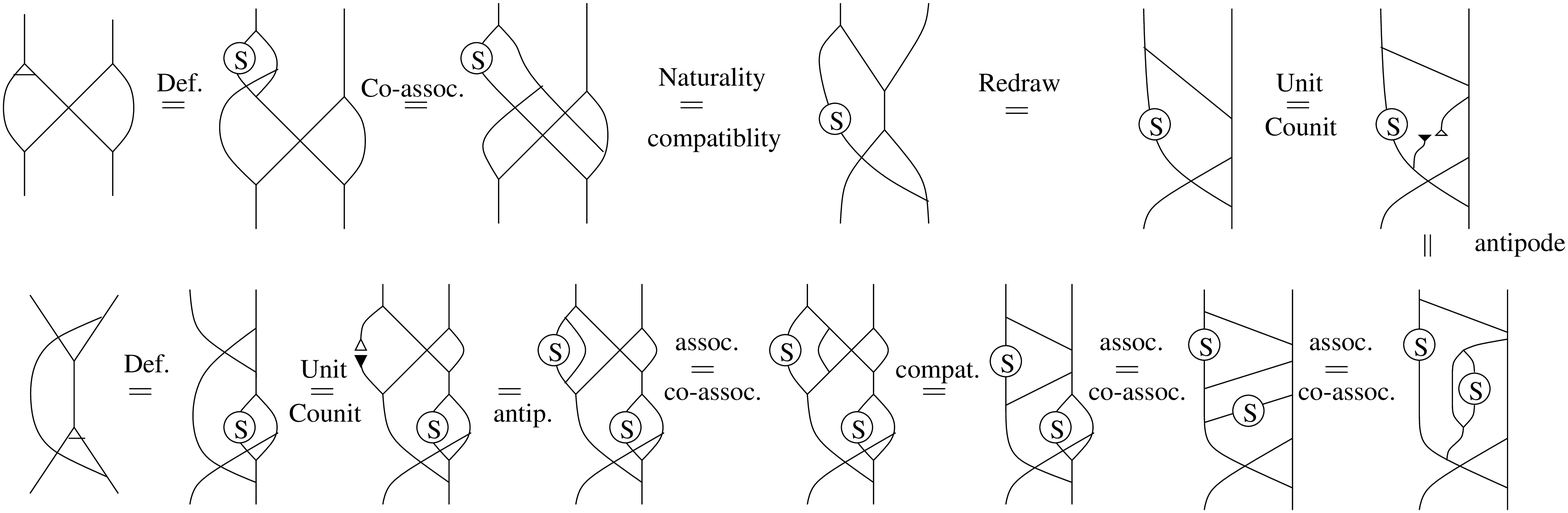}
}
\end{center}
\caption{Adjoint condition $(2)$}
\label{hopfYBElem}
\end{figure}

It is known that these properties are satisfied, and proofs are found
in \cite{Henn,Woro}.
Here we include diagrammatic proofs for reader's convenience
in Fig.~\ref{HopfQlemma} and Fig.~\ref{hopfYBElem}, respectively.
$\Box$

\begin{definition}{\rm
We call the above equalities (\ref{ad1eq}) and (\ref{ad2eq}) the
{\it adjoint conditions}.
} \end{definition}

\begin{remark}{\rm
The equality (\ref{ad1eq}) is equivalent to the fact that the
adjoint map defines an algebra action of $H$ on itself (see
\cite{MajidGreen}). Specifically, $(a \lt
b)\lt c = a \lt (bc)$ for any $a,b,c \in H$, where $\lt$ denotes
the right action defined by the adjoint: $a \lt b=\ad(a \otimes
b)$.
The equality (\ref{ad2eq}) can be similarly rewritten as:
$$a_{(1)} \lt b_{(1)} \otimes a_{(2)}  b_{(2)}=
( a \lt b_{(2)} )_{(1)} \otimes b_{(1)} ( a \lt b_{(2)} )_{(2)} . $$
} \end{remark}

\begin{remark}{\rm
It was pointed out to us by
Sommerhaeuser
that the induced $R$-matrix $R_{\ad}$ is invertible with inverse
$$R^{-1}_\ad (b\otimes a) =  b_{(3)} a S^{-1}(b_{(2)}) \otimes b_{(1)} .$$
} \end{remark}

\section{Deformations of  the Adjoint Map}\label{deformsec}

We follow the exposition in \cite{MrSt} for deformation of
bialgebras to propose  a similar deformation theory for the
adjoint map. In light of Lemma~\ref{adYBElem}, we deform the two
equalities (\ref{ad1eq}) and (\ref{ad2eq}).
Let $H$ be a Hopf algebra and $\ad$ its adjoint map.

\begin{definition} {\rm
A deformation of $(H, \ad) $ is a pair $(H_t, \ad_t)$ where $H_t$
is a $k[[t]]$-Hopf algebra
given by $H_t=H \otimes k[[ t ]]$ with all Hopf algebra structures
inherited by extending those on $H_t$
with the identity on the
$k[[t]]$ factor (the trivial deformation as a Hopf algebra), with
a deformations of $\ad$ given by $\ad_t= \ad + t \ad_1 + \cdots +
t^n \ad_n + \cdots : H_t \otimes H_t \rightarrow H_t$ where
$\ad_i: H \otimes H \rightarrow H $,
 $i=1, 2, \cdots$, are
maps. } \end{definition}

Suppose $\bar{\ad}=\ad + \cdots + t^n \ad_n$
 satisfies  the adjoint
 conditions (equalities  (\ref{ad1eq}) and (\ref{ad2eq})) mod $t^{n+1}$,
and suppose
that  there exist $\ad_{n+1}: H \otimes H \rightarrow H$  such that
$\bar{\ad}+t^{n+1} \ad_{n+1}$
satisfies the adjoint  conditions mod $t^{n+2}$.
Define  $\xi_1 \in \Hom(H^{\otimes 3}, H)$ and
$\xi_2\in \Hom(H^{\otimes 2}, H^{\otimes 2})$
by
\begin{eqnarray*}
\bar{\ad} (\bar{\ad} \otimes 1) - \bar{\ad}(1 \otimes \mu)
&=& t^{n+1} \xi_1 \quad {\rm mod}\  t^{n+2} , \label{ad2d1} \\
 (\bar{\ad} \otimes \mu)(1 \otimes \tau \otimes 1)(\Delta
\otimes \Delta) \hspace{2in} \\
 - (1 \otimes \mu)(\tau \otimes 1)(1
\otimes \Delta)(1 \otimes \bar{\ad})(\tau \otimes 1)(1 \otimes
\Delta)&=&  t^{n+1} \xi_2 \quad {\rm mod}\  t^{n+2} .
 \label{ad2d2}
\end{eqnarray*}
For the first adjoint condition (\ref{ad1eq})
 of $\bar{\ad}+t^{n+1} \ad_{n+1}$ mod $t^{n+2}$ we obtain:
$$( \bar{\ad}+t^{n+1} \ad_{n+1})
((\bar{\ad}+t^{n+1} \ad_{n+1}) \otimes 1) - (\bar{\ad}+t^{n+1}
\ad_{n+1})(1 \otimes \mu) =0 \  {\rm mod}\  t^{n+2} $$
which is equivalent by degree calculations to:
$$
{\ad} (\ad_{n+1} \otimes 1)
+ \ad_{n+1} ({\ad} \otimes 1)
 - \ad_{n+1}(1 \otimes \mu)
=\xi_1. $$
For the second adjoint condition (\ref{ad2eq})
 of $\bar{\ad}+t^{n+1} \ad_{n+1}$ mod $t^{n+2}$ we obtain:
\begin{eqnarray*}
\lefteqn{ ( (\bar{\ad}+t^{n+1} \ad_{n+1})
 \otimes \mu)(1 \otimes \tau \otimes 1)(\Delta
\otimes \Delta)} \\
& - & (1 \otimes \mu)(\tau \otimes 1)(1 \otimes \Delta)(1 \otimes
(\bar{\ad}+t^{n+1} \ad_{n+1} ))(\tau \otimes
1)(1 \otimes \Delta) \\
&=&  0 \quad  {\rm mod}\  t^{n+2}
\end{eqnarray*}
which is equivalent by degree calculations to:
\begin{eqnarray*}
\lefteqn{ ( \ad_{n+1}
 \otimes \mu)(1 \otimes \tau \otimes 1)(\Delta
\otimes \Delta)} \\
& - & (1 \otimes \mu)(\tau \otimes 1)(1 \otimes \Delta)(1 \otimes
\ad_{n+1} )(\tau \otimes 1)(1 \otimes \Delta) \quad = \quad \xi_2
.
\end{eqnarray*}

In summary we proved the following:
\begin{lemma} \label{deformlem}
The map   $\bar{\ad}+t^{n+1} \ad_{n+1}$  satisfies the
adjoint conditions mod $t^{n+2}$ if and only if
\begin{eqnarray*}
{\ad} (\ad_{n+1} \otimes 1)
+ \ad_{n+1} (
{\ad} \otimes 1)
 - \ad_{n+1}(1 \otimes \mu) &=& \xi_1,  \\
{\rm and } \quad  ( \ad_{n+1}
 \otimes \mu)(1 \otimes \tau \otimes 1)(\Delta
\otimes \Delta) \hspace{1in} & &  \\
 -  (1 \otimes \mu)(\tau \otimes 1)(1
\otimes \Delta)(1 \otimes  \ad_{n+1} )(\tau \otimes 1)(1 \otimes
\Delta) \quad &=&   \xi_2 .
\end{eqnarray*}

\end{lemma}

\section{Differentials and Cohomology}
\label{Diff}

\subsection{Chain Groups}

We define chain groups, for positive integers $n$, $n>1$, and $i
=1, \ldots, n$  by:
\begin{eqnarray*}
C^{n, i}_{\rm ad} (H;H) & =&  \Hom(H^{\otimes (n+1-i)}, H^{\otimes i} ), \\
C^n_{\rm ad} (H;H) & =& \displaystyle{\oplus_{ i>0,\  i \leq n+1-i}
\ C^{n, i}_{ad} (H;H)}.
\end{eqnarray*}
Specifically, chain groups in low dimensions of our concern are:
\begin{eqnarray*}
C^2_{\rm ad}(H;H)&=& \Hom(H^{\otimes 2}, H), \\
C^3_{\rm ad}(H;H)&=& \Hom(H^{\otimes 3}, H)\oplus
\Hom(H^{\otimes 2}, H^{\otimes 2}).
\end{eqnarray*}

For $n=1$, define
$$C^1_{\rm ad}(H;H)=  \{ f \in \Hom_k (H,H) \ | \
f \mu = \mu (f \otimes 1) + \mu (1 \otimes f), \
\Delta f = (f \otimes 1 ) \Delta + (1 \otimes f) \Delta \ \}.  $$

In the remaining sections we will define differentials that are
homomorphisms between the chain groups:
$$d^{n, i} :  C^n_{\rm ad} (H;H)
 \rightarrow C^{n+1, i}_{\rm ad}(H;H) (= \Hom(H^{\otimes (n+2-i)}, H^{\otimes i} ) )$$
that will be  defined individually for $n=1,2,3$ and for $i$ with $2i \leq n+1$,
and
\begin{eqnarray*}
 D_1 & =& d^{1, 1}: C^1_{\rm ad} (H;H) \rightarrow C^2_{\rm ad} (H;H),\\
 D_2&=& d^{2, 1} + d^{2, 2}:  C^2_{\rm ad} (H;H) \rightarrow C^3_{\rm ad} (H;H),\\
 D_3&=& d^{3, 1} +d^{3, 2} +d^{3, 3}:
  C^3_{\rm ad} (H;H) \rightarrow C^3_{\rm ad} (H;H).
\end{eqnarray*}

\subsection{First Differentials}

By analogy with the differential for multiplication, we
make the following definition:

\begin{figure}[htb]
\begin{center}
\mbox{
\epsfxsize=2.5in
\epsfbox{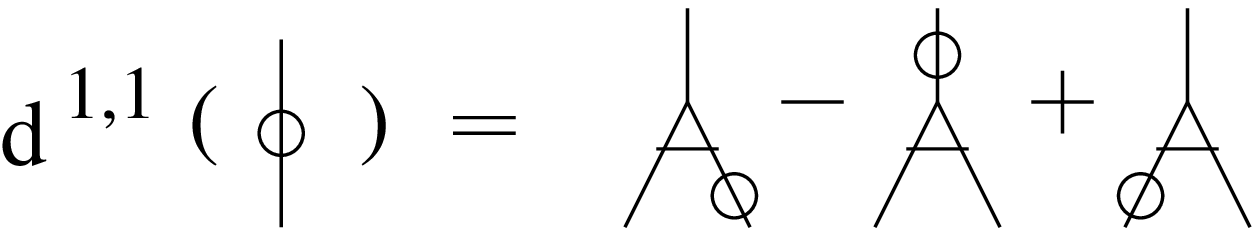}
}
\end{center}
\caption{The $1$-differential }
\label{d11}
\end{figure}

\begin{definition}{\rm
The first differential
$$d^{1,1}:
C^{1}
_{\rm ad} (H;H)
\rightarrow
C^{2,1}_{\rm ad} (H;H)
$$
is defined by
$$d^{1,1}(f)= \ad (1 \otimes f) - f\ad + \ad(f \otimes 1). $$
} \end{definition}
Diagrammatically, we represent $d^{1,1}$ as depicted in
Fig.~\ref{d11},
where a $1$-cochain is represented by a circle on a string.

\begin{figure}[htb]
\begin{center}
\mbox{
\epsfxsize=.4in 
\epsfbox{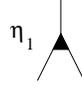}
}
\end{center}
\caption{A diagram for a $2$-cochain }
\label{co2}
\end{figure}

\begin{figure}[htb]
\begin{center}
\mbox{
\epsfxsize=3in
\epsfbox{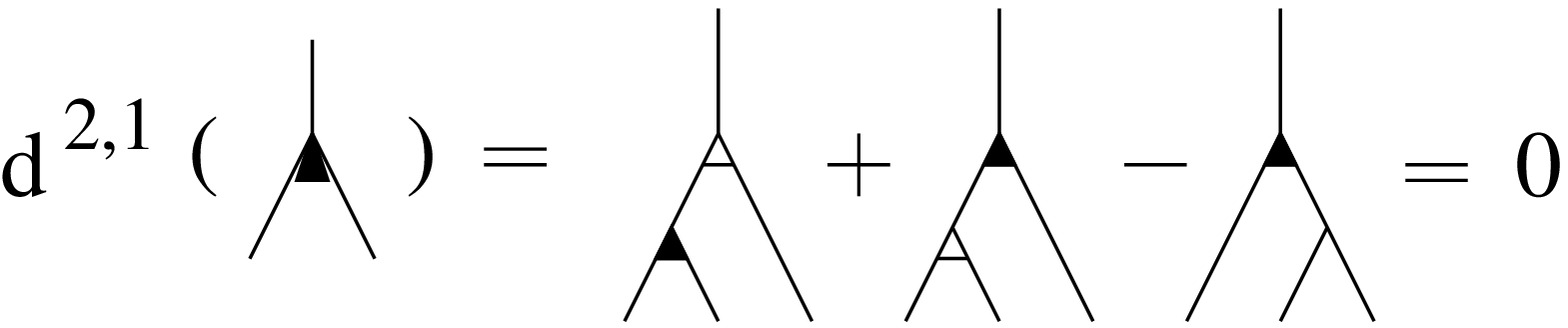}
}
\end{center}
\caption{The $2$-cocycle condition, Part I }
\label{d21}
\end{figure}

\begin{figure}[htb]
\begin{center}
\mbox{
\epsfxsize=2.5in
\epsfbox{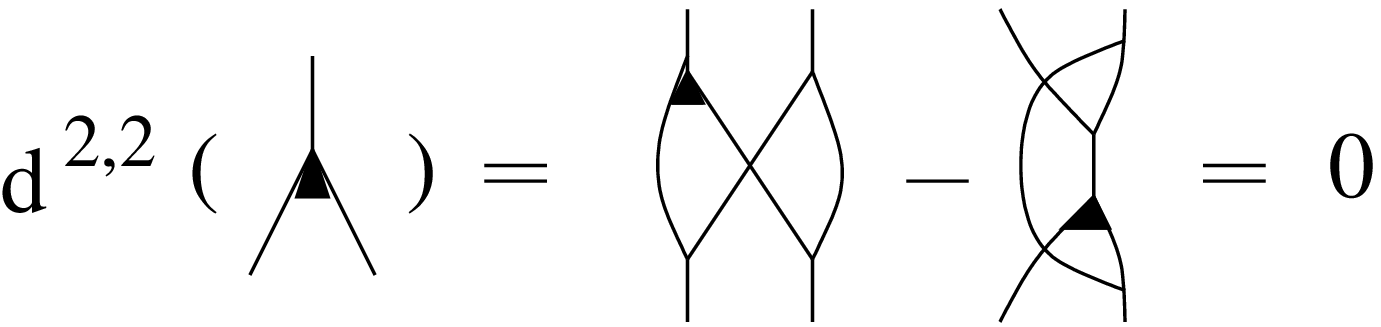}
}
\end{center}
\caption{The  $2$-cocycle condition, Part II}
\label{d22}
\end{figure}

\subsection{Second Differentials}

\begin{definition}{\rm
Define the second  differentials by:
\begin{eqnarray*}
d^{2,1}_\ad (\phi) &=& \ad (\phi \otimes 1) + \phi (\ad \otimes 1) - \phi(1 \otimes \mu) , \\
d^{2,2}_\ad (\phi) &=& (\phi \otimes \mu)(1 \otimes \tau \otimes
1)(\Delta \otimes \Delta) - (1 \otimes \mu)(\tau \otimes 1)(1
\otimes \Delta)(1 \otimes \phi) (\tau \otimes 1) (1 \otimes
\Delta).
\end{eqnarray*}
} \end{definition}
Diagrams for $2$-cochain and $2$-differentials are  depicted in Fig.~\ref{co2},
\ref{d21}, and \ref{d22}, respectively.

\begin{figure}[htb]
\begin{center}
\mbox{
\epsfxsize=4in
\epsfbox{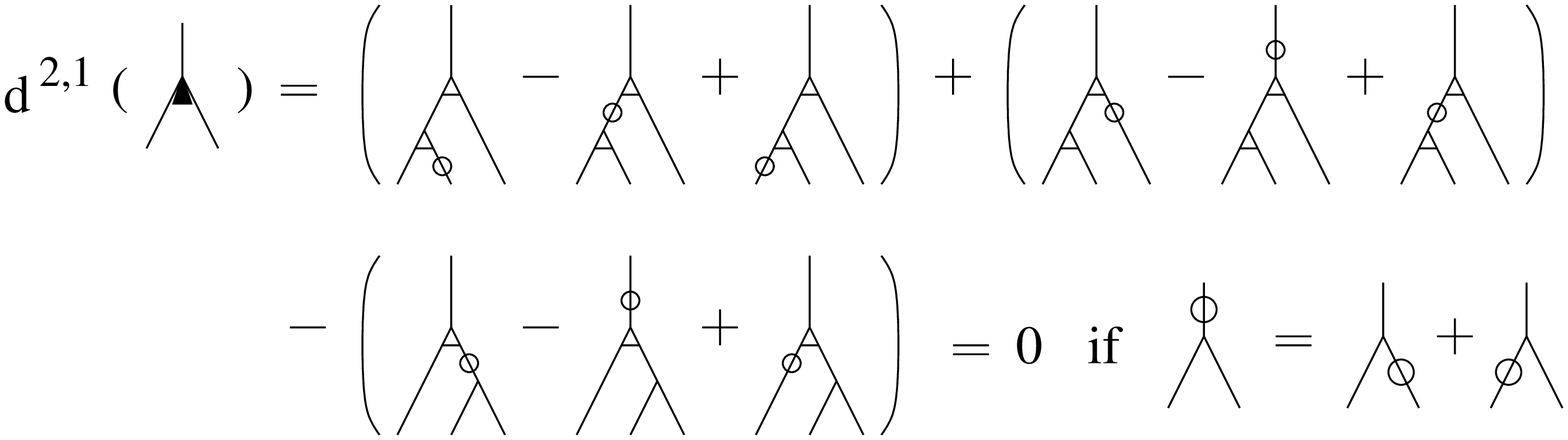}
}
\end{center}
\caption{The $2$-cocycle condition for a $2$-coboundary, Part I }
\label{dd21}
\end{figure}

\begin{figure}[htb]
\begin{center}
\mbox{
\epsfxsize=5.5in
\epsfbox{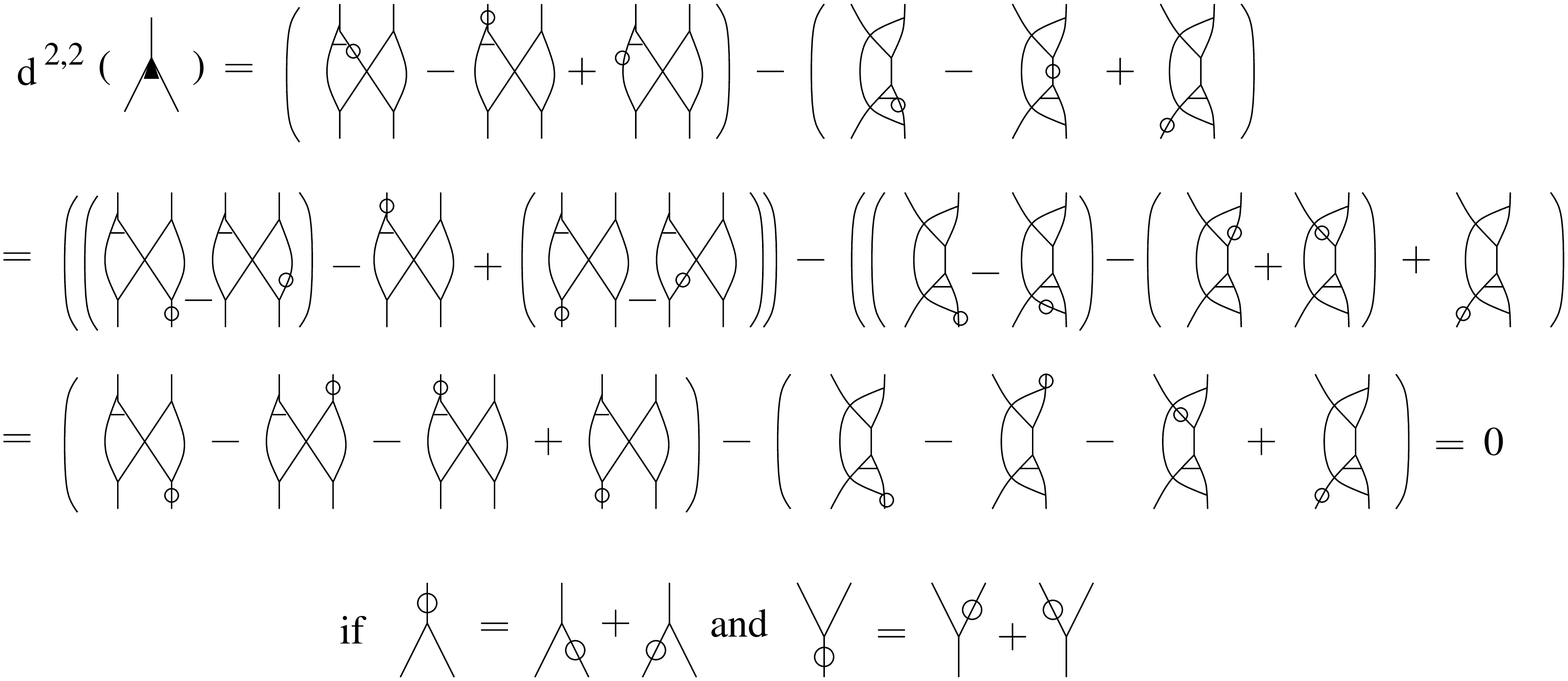}
}
\end{center}
\caption{The  $2$-cocycle condition for a $2$-coboundary, Part II}
\label{dd22}
\end{figure}

\begin{theorem}\label{dd2thm}
$D_2 D_1=0$.
\end{theorem}
{\it Proof.\/} This follows from
direct calculations, and
can be seen from diagrams in Figs.~\ref{dd21} and  \ref{dd22}.
$\Box$

\begin{figure}[htb]
\begin{center}
\mbox{
\epsfxsize=2.5in
\epsfbox{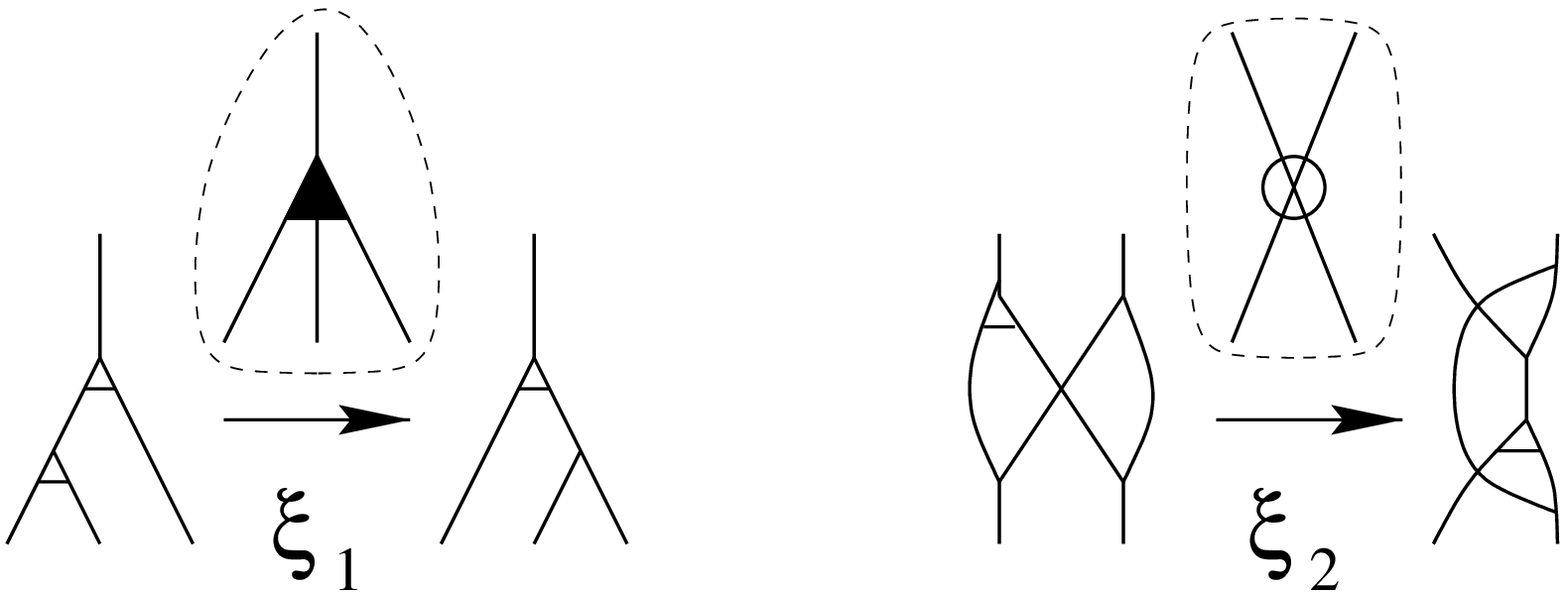}
}
\end{center}
\caption{ Diagrams for $3$-cochains}
\label{co3}
\end{figure}

\begin{figure}[htb]
\begin{center}
\mbox{
\epsfxsize=2.2in
\epsfbox{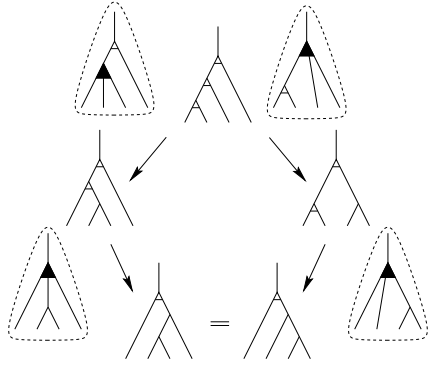}
}
\end{center}
\caption{The first $3$-differential $d^{3,1}$} \label{d31}
\end{figure}

\begin{figure}[htb]
\begin{center}
\mbox{
\epsfxsize=3.5in
\epsfbox{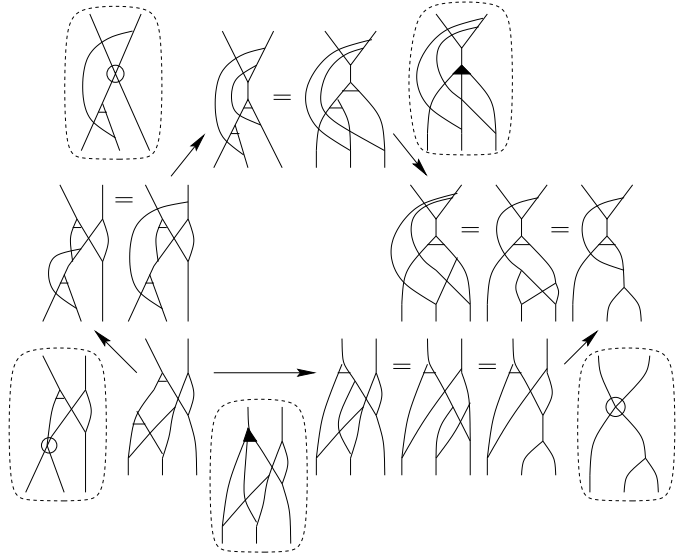}
}
\end{center}
\caption{The second $3$-differential $d^{3,2}$} \label{d32}
\end{figure}

\begin{figure}[htb]
\begin{center}
\mbox{
\epsfxsize=2.5in
\epsfbox{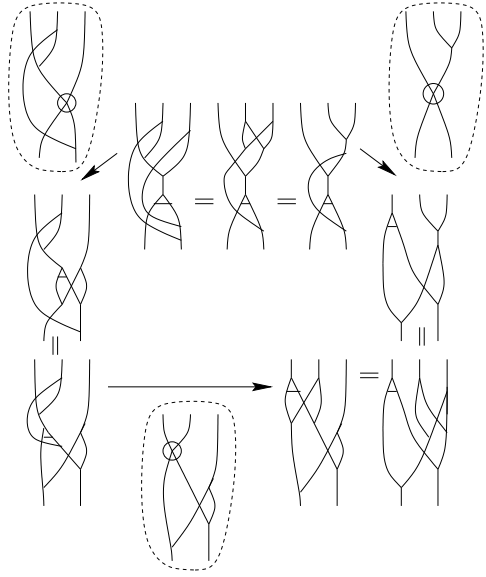}
}
\end{center}
\caption{ The third $3$-differential $d^{3,3}$} \label{d33}
\end{figure}

\subsection{Third Differentials}

\begin{definition}\label{third}{\rm
We define $3$-differentials as follows.
Let $\xi_i \in C^{3,i}(H;H)$ for $i=1,2$.
Then
\begin{eqnarray*}
d^{3,1}_\ad(\xi_1, \xi_2)&=& \ad (\xi_1 \otimes 1) +  \xi_1(1
\otimes \mu \otimes 1)
- \xi_1( \ad  \otimes 1^2 + 1^2 \otimes \mu ) , \\
d^{3,2}_\ad(\xi_1, \xi_2)&=& ( \ad \otimes \mu)(1 \otimes \tau
\otimes 1)(1^2 \otimes \Delta)(\xi_2 \otimes 1)
+(1 \otimes \mu)(\tau \otimes 1) (1 \otimes \xi_2) (R_\ad \otimes 1) \\
& & + (1 \otimes \mu)(1^2 \otimes \mu)(\tau \otimes 1^2)(1 \otimes
\tau \otimes 1)
(1^2 \otimes \Delta)
((1^2 \otimes \xi_1) \\
& & \quad \cdot (1 \otimes \tau \otimes 1^2)(\tau \otimes 1^3)(1^2
\otimes \tau \otimes 1) (1 \otimes \Delta \otimes \Delta) ) \\
& & - (\xi_1 \otimes \mu) (1^{2} \otimes \tau \otimes 1)(1^{2}\otimes \mu \otimes 1^{2})(1 \otimes \tau \otimes 1^{3})(\Delta \otimes \Delta \otimes \Delta)
 - \xi_2 (1 \otimes \mu),\\
 d^{3,3}_\ad(\xi_1, \xi_2)&=&
(1 \otimes \mu \otimes 1) (\tau  \otimes 1^2)(1 \otimes \Delta
\otimes 1) (1 \otimes \xi_2 )
(\tau \otimes 1)
(1 \otimes \Delta)
\\ & & + (\xi_2 \otimes \mu) (1 \otimes \tau \otimes 1)(\Delta \otimes \Delta)
- (1 \otimes \Delta) \xi_2 .
\end{eqnarray*}
} \end{definition}
Diagrams for  $3$-cochains are depicted in Fig.~\ref{co3}. See
Figs.~\ref{d31}, \ref{d32}, and \ref{d33}  for the diagrammatics
for $d^{3,1}$, $d^{3,2}$ and $d^{3,3}$, respectively.

\begin{figure}[htb]
\begin{center}
\mbox{
\epsfxsize=3.2in
\epsfbox{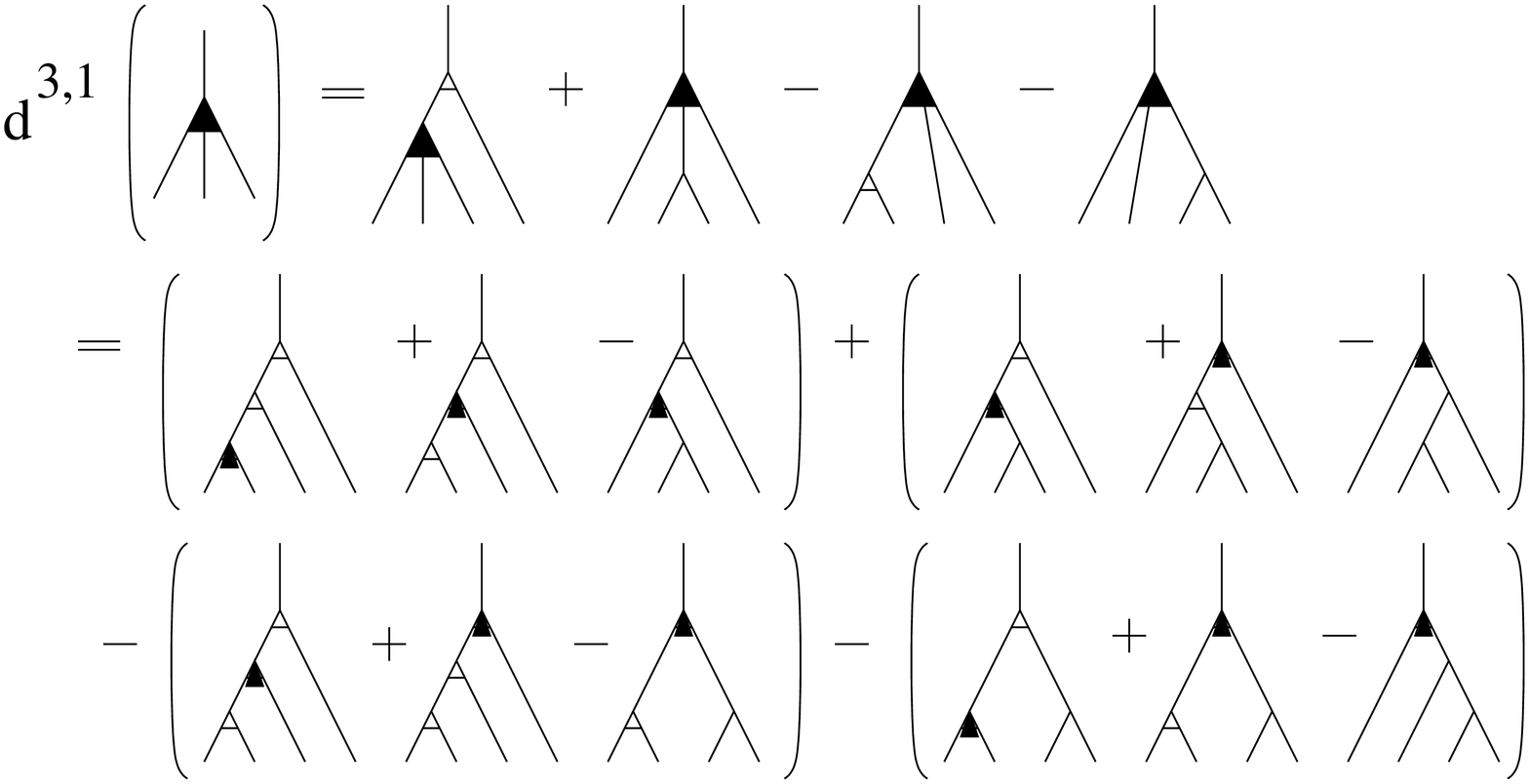}
}
\end{center}
\caption{$d^{3,1}(d^{2,1}) = 0$}
\label{dd31A}
\end{figure}

\begin{figure}[htb]
\begin{center}
\mbox{
\epsfxsize=4.5in
\epsfbox{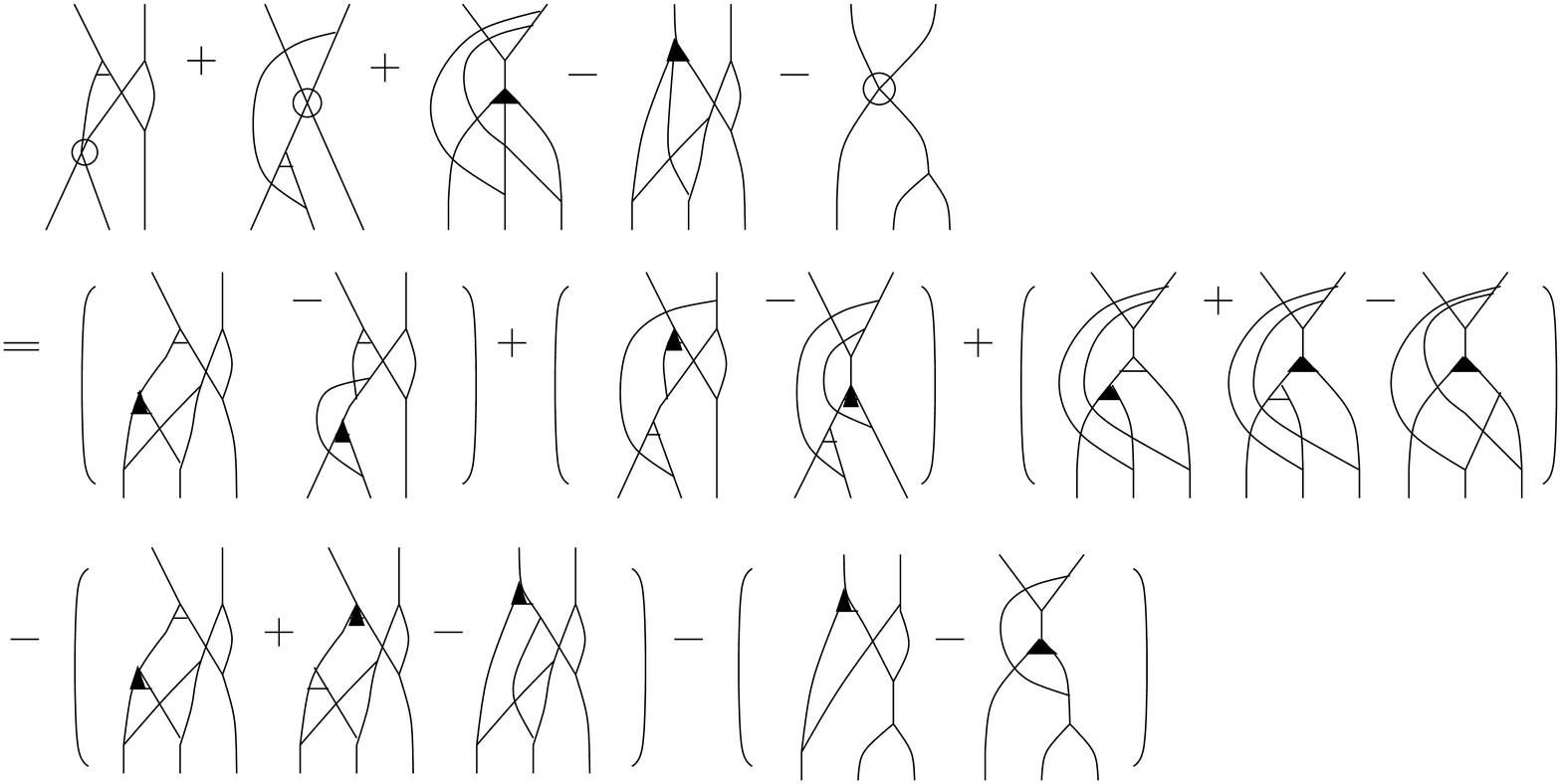}
}
\end{center}
\caption{$d^{3,2}(d^{2,1}, d^{2,2}) = 0$}
\label{dd32A}
\end{figure}

\begin{figure}[htb]
\begin{center}
\mbox{
\epsfxsize=4in
\epsfbox{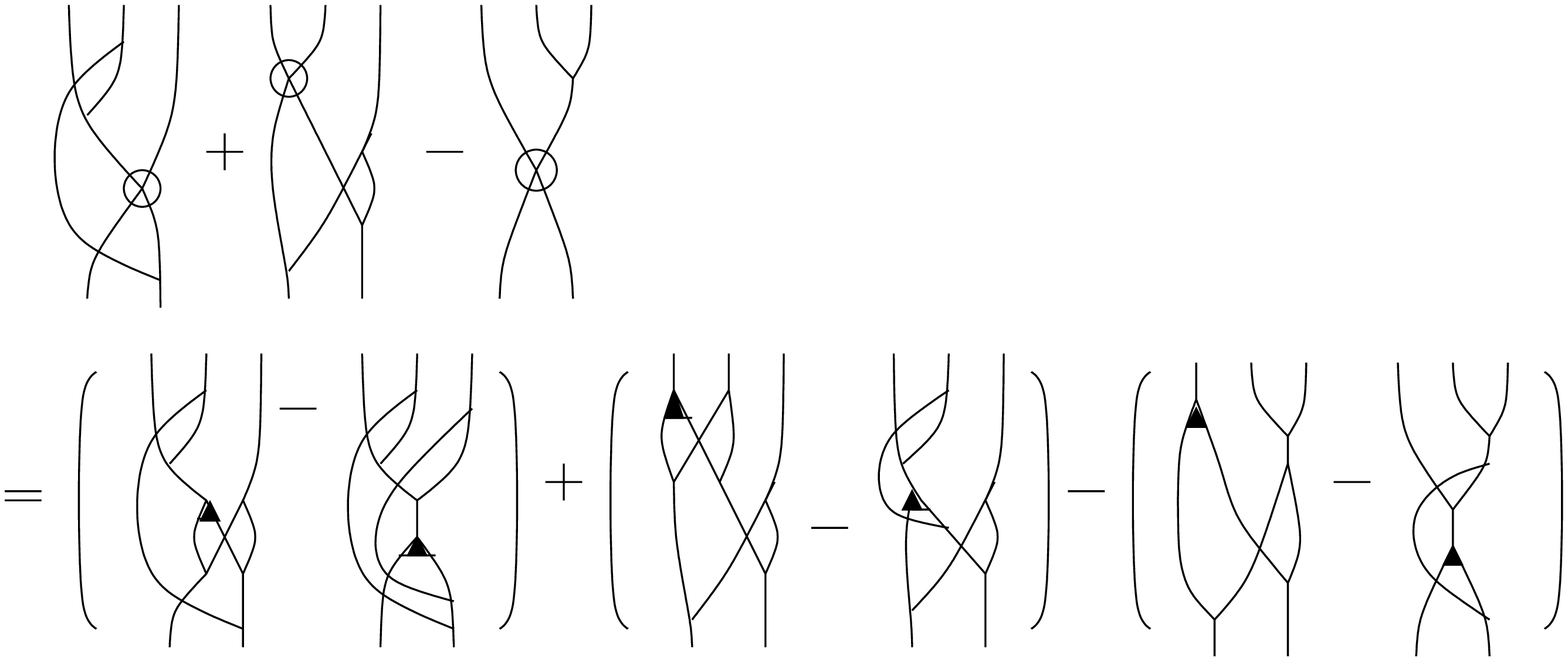}
}
\end{center}
\caption{$d^{3,3}(d^{2,2}) = 0$}
\label{dd33A}
\end{figure}

\begin{theorem}\label{dd3thm}
$D_3 D_2 = 0 $.
\end{theorem}
{\it Proof.\/} The proof follows from
direct calculations that are indicated in
 Figs.~\ref{dd31A}, \ref{dd32A} and  \ref{dd33A}.
 We demonstrate how to recover algebraic calculations from
 these diagrams for
 the part $(d^{3,3} d^{2,2})(\eta_1)=0$
  for any $\eta_1\in C^2(H;H)$.
  This is indicated in Fig.~\ref{dd33A}, where subscripts $\ad$ are
 suppressed for simplicity.
 Let
 $\xi_2=d^{2,2}(\eta_1)\in C^{3,2}(H;H)$
 (note that $\xi_1=d^{2,1}(\eta_1)\in C^{3,1}(H;H)$ does not land in the domain of
  $d^{3,3}$).
 The first line of Fig.~\ref{dd33A} represents the definition of the differential
\begin{eqnarray*}
 d^{3,3}_\ad(\xi_1, \xi_2)&=&
(1 \otimes \mu \otimes 1) (\tau  \otimes 1^2)(1 \otimes \Delta
\otimes 1) (1 \otimes \xi_2 ) (\tau \otimes 1) (1 \otimes \Delta)
\\ & & + (\xi_2 \otimes \mu) (1 \otimes \tau \otimes 1)(\Delta \otimes \Delta)
- (1 \otimes \Delta) \xi_2
 \end{eqnarray*}
 where each term represents
 each connected diagram.
 The first parenthesis of the second line represents
that
 \begin{eqnarray*}
\xi_2=d^{2,2} (\eta_1) &=& (\eta_1 \otimes \mu)(1 \otimes \tau
\otimes 1)(\Delta \otimes \Delta) - (1 \otimes \mu)(\tau \otimes
1)(1 \otimes \Delta)(1 \otimes
\eta_1
) (\tau \otimes 1) (1 \otimes
\Delta)
\end{eqnarray*}
is substituted in the first term
$$(1 \otimes \mu \otimes 1)
(\tau  \otimes 1^2)(1 \otimes \Delta \otimes 1)
(1 \otimes \xi_2 )
(\tau \otimes 1)
(1 \otimes \Delta) .$$
When these two maps are applied to a general element $x \otimes y \in H \otimes H$,
the results are computed as
 \begin{eqnarray*}
& & \eta_1( x_{(1)} \otimes y_{(2)(1)} )_{(1)}  \otimes
 y_{(1)}  \eta_1( x_{(1)} \otimes y_{(2)(1)} )_{(2)} \otimes
 x_{(2)} y_{(2)(2)}, \\
 & & -
 \eta_1 (x \otimes y_{(2)(2)} )_{(1)(1)} \otimes
 y_{(1)}\eta_1 (x \otimes y_{(2)(2)} )_{(1)(2)} \otimes
 y_{(2)(1)}\eta_1 (x \otimes y_{(2)(2)} )_{(2)}.
 \end{eqnarray*}
By coassociativity applied to $y$ and $\eta_1 (x \otimes y_{(2)(2)} )$,
the
second 
 term is equal to
$$\eta_1 (x \otimes y_{(2)} )_{(1)} \otimes
 y_{(1)(1)}\eta_1 (x \otimes y_{(2)} )_{(2)(1)} \otimes
 y_{(1)(2)}\eta_1 (x \otimes y_{(2)} )_{(2)(2)}, $$
 which is equal, by compatibility, to
 $$\eta_1 (x \otimes y_{(2)} )_{(1)} \otimes
( y_{(1)}\eta_1 (x \otimes y_{(2)} )_{(2)} )_{(1)} \otimes
 (y_{(1)}\eta_1 (x \otimes y_{(2)} )_{(2)})_{(2)} . $$
 This last term is represented exactly
 by  the
  last term in the second line of
 Fig.~\ref{dd33A}, and therefore is cancelled.
 The map represented by the second term in the second line of  Fig.~\ref{dd33A}
 cancels with the third term by coassociativity,
 and the fourth term cancels with the sixth by coassociativity applied twice
 and
compatibility once.
 Other cases (Figs.~\ref{dd31A}, \ref{dd32A}) are computed similarly.
$\Box$

\subsection{Cohomology Groups}

For convenience define
$C^0(H;H)=0$, $D_0=0: C^0(H;H)\rightarrow C^1(H;H)$.

Then Theorems~\ref{dd2thm}, \ref{dd3thm} are summarized as:
\begin{theorem}
${\cal C}=(C^n, D_n)_{n=0,1,2,3}$ is a chain complex.
\end{theorem}

This enables us to define:
\begin{definition}{\rm
The {\it adjoint } $n$-coboundary, cocycle and cohomology groups
are defined by:
\begin{eqnarray*}
B^{n}(H;H) &=& {\rm Image}(D_{n-1}) , \\
Z^{n}(H;H) &=& {\rm Ker}(D_n) , \\
H^{n}(H;H) &=& Z^{n}(H;H) / B^{n}(H;H)
\end{eqnarray*}
for $n=1,2,3$.
} \end{definition}

\section{Examples} \label{Examples}

\subsection{Group Algebras}\label{gpalgsubsec}

Let $G$ be a group and $H=kG$
be its
group algebra with
the coefficient field $k$
(char $k \ne 2$).
Then $H$ has a Hopf algebra structure induced
from the group operation as multiplication, $\Delta(x)=x \otimes x$ for
basis elements $x \in G$, and the antipode induced from $S(x)=x^{-1}$
for $x \in G$.
Here and below, we denote the conjugation
action on a group $G$ by $x \lt y:= y^{-1} x y$.
Note that this defines a quandle structure on $G$; see \cite{Joyce}.

\begin{lemma}\label{kGC1lem}
$C^1_{\rm ad}(kG;kG)=0$.
\end{lemma}
{\it Proof.\/} For any given $w \in G$ write
$\displaystyle{f(w)=\sum_{u \in G} a_u(w) u}$, where $a: G
\rightarrow k$ is a function. Recall the definition
$$C^1_{\rm ad}(H;H)=  \{ f \in \Hom_k (H,H) \ | \
f \mu = \mu (f \otimes 1) + \mu (1 \otimes f), \
\Delta f = (f \otimes 1 ) \Delta + (1 \otimes f) \Delta \ \}. $$
The LHS of the second condition is written as
$$\Delta f(w)=\Delta (\sum_{u} a_u(w) u)=\sum_{u} a_u(w) u \otimes u $$
and the RHS is written as
$$ ((f \otimes 1 ) \Delta + (1 \otimes f) \Delta )(w)=
( (f \otimes 1 )+ (1 \otimes f))(w \otimes w) = \sum_h a_h(w) (h
\otimes w) + \sum_v a_v (w)(w \otimes v).$$
For a given $w$, fix $u$ and then compare the coefficients of $u
\otimes u$.
  In the LHS we have $a_u(w)$, 
while on the RHS
$w=u$, and furthermore
$w=h=v$ for $u \otimes u$.  Thus the  diagonal
coefficient must satisfy  $a_w
(w)=a_w(w)+a_w(w)$, so that $a_w(w)=0$ since char $k\ne2$.
In the case $w \neq u$, neither term of $h \otimes w$ nor $w
\otimes v$ is equal to $u \otimes u$, hence $a_u(w)=0$. $\Box$

\begin{lemma}\label{G2lem}
For $x, y \in G$, write $\displaystyle{\phi(x, y)=\sum_u a_u(x, y)
u}$, where $a: G \times G \rightarrow k$. Then the induced linear
map $\phi : kG \otimes kG \rightarrow kG$ is in $Z^2_\ad (kG; kG)$
if and only if $a$ satisfies
$$ a_{x\lt y} (x, y) + a_{(x \lt y)\lt z}(x \lhd y, z) - a_{(x \lt y)\lt z}(x, yz)=0 $$
for any $x,y,z \in G$.
\end{lemma}
{\it Proof.\/}
The first $2$-cocycle condition for $\phi : kG
\otimes kG \rightarrow kG$ is written by:
$$z^{-1} \phi (x  \otimes  y)  z
+ \phi(y^{-1}xy  \otimes  z) - \phi(x \otimes  yz)=0$$ for basis
elements $x, y, z \in G$.
The  second is formulated by
$$ {\rm LHS}=\phi(x \otimes  y)\otimes xy= \sum_u a_u(x, y) (u\otimes xy ),
\quad {\rm RHS}=\sum_w a_w(x, y) (w\otimes y w ) . $$ 
 They have the common term $u \otimes xy$ for
$w=y^{-1}xy=u$, and otherwise they are different terms. Thus we
obtain $a_w(x,y)=0$ unless $w=y^{-1}xy$. For these terms, the
first condition becomes
$$ z^{-1} (a_{y^{-1}xy} (x, y) y^{-1} x y ) z
+ a_{z^{-1}y^{-1}xyz}(y^{-1}xy, z) z^{-1}y^{-1}xyz -
a_{z^{-1}y^{-1}xyz}(x, yz)z^{-1}y^{-1}xyz=0 $$ and the result
follows. $\Box$

\begin{remark}\label{G2rem}{\rm
In the preceding proof,
since the term $a_w(x, y)=0$ unless $w=x \lt y$,
let $a_{x \lt y}(x,y)=a(x,y)$.
Then the condition stated becomes
$$a(x, y)+a(x \lt y, z)-a(x, yz)=0.$$
}\end{remark}

\begin{proposition}\label{Gp3cocyProp}
Let $G$ be a group.
Let $(\xi_1, \xi_2) \in C^3_{\ad}(kG; kG)$, where $\xi_1$ is the
map that is defined by linearly extending $\xi_1(x \otimes y
\otimes z) = \displaystyle{\sum_{u \in G} c_{u} (x,y,z) u }$. Then
$(\xi_1 , \xi_2) \in Z^3_{\ad}(kG;kG)$ if and only if $\xi_2 = 0$
and the coefficients satisfy the following properties:
\begin{itemize}
\item[(a)] $c_u(x,y,z) = 0$ if $u \neq z^{-1}y^{-1}xyz$ and
\item[(b)] $c(x,y,z)=c_{ z^{-1}y^{-1}xyz} (x,y,z) $ satisfies
$$c(x, y, z)+ c(x, yz, w)=c(y^{-1} x y,
z, w) + c(x, y, zw).$$
\end{itemize}
\end{proposition}
{\it Proof.\/} Suppose $(\xi_1, \xi_2) \in Z^3_{\ad}(kG;kG)$. Let
$\xi_2$ be the map that is defined by linearly extending
$\displaystyle{\xi_2(x \otimes y)= \sum_{u, v \in G} a_{u,v} (x,y)
u \otimes v }$. Then the third $3$-cocycle condition from
Definition \ref{third} gives: (abbreviating
$a_{u,v}(x,y)=a_{u,v}$)
\begin{eqnarray*}
\lefteqn{d^{3,3}(\xi_1, \xi_2)(x \otimes y)}\\
&=&
\sum_{u_1, v_1}  a_{u_1,v_1}  (u_1 \otimes y  u_1 \otimes v_1)
+ \sum_{u_2, v_2}  a_{u_2,v_2} (u_2 \otimes v_2 \otimes xy )
- \sum_{u_3, v_3}  a_{u_3,v_3} (u_3 \otimes v_3 \otimes v_3  ) \quad = 0 .
\end{eqnarray*}
We first consider terms in which the third tensorand is $xy$. From
the third summand, this forces the second tensorand to be $xy$, so
we collect the terms of the form $(u \otimes xy \otimes xy)$. This
gives:
$$\sum _{u} (a_{u, xy} + a_{u, xy} - a_{u, xy})(u \otimes xy
\otimes xy) = 0,$$ which implies $a_{u, xy} = 0$ for all $u \in
G$.  The remaining terms are $$ \sum_{u_1, v_1\neq xy}
a_{u_1,v_1}  (u_1 \otimes y  u_1 \otimes v_1) + \sum_{u_2, v_2\neq
xy}  a_{u_2,v_2} (u_2 \otimes v_2 \otimes xy ) - \sum_{u_3,
v_3\neq xy}  a_{u_3,v_3} (u_3 \otimes v_3 \otimes v_3  )  = 0 .
$${}From the second
sum
we obtain $a_{u, v}(x, y)=0$ for $v\neq
xy$. In conclusion, if $d^{3,3}(\xi_1, \xi_2) = 0$ for $kG$ then
$\xi_2 = 0$.

We now consider $d^{3,2}(\xi_1, \xi_2)$, with $\xi_2 = 0$.  Let
$\xi_1$ be the map that is defined by linearly extending
$\displaystyle{\xi_1(x \otimes y \otimes z)=\sum_{u \in G} c_{u}
(x,y,z) u }$ for $x,y,z \in G$. The second $3$-cocycle condition
from Definition \ref{third}, with $\xi_2=0$, is
$\displaystyle{\sum_{u } c_{u} u \otimes yz u =\sum_{v } c_{v} v
\otimes xyz}$.  In order to combine like terms, we need $yzu =
xyz,$ meaning $u = z^{-1} y^{-1} xyz$. Thus, $c_u(x , y, z)=0$
except in the case when $u = z^{-1}y^{-1}xyz$. In this case, we
obtain $\xi_1(x \otimes y \otimes z)= c(x,y,z) z^{-1}y^{-1}xyz
\otimes xyz$ where $c(x,y,z)=c_{z^{-1}y^{-1}xyz} (x,y,z)$.

Finally we consider the first $3$-cocycle condition from
Definition \ref{third}, which is formulated for basis elements by
$$w^{-1} \ \xi_1(x \otimes y \otimes z)\  w +  \xi_1(x\otimes yz\otimes w)
=  \xi_1(x \lt y \otimes z\otimes w)+ \xi_1(x\otimes y\otimes zw).
$$ Substituting in the formula for $c(x,y,z)$ which we found above, we obtain
$$c(x, y, z)+ c(x, yz, w)=c(y^{-1} x y,
z, w) + c(x, y, zw).$$
This is a group $3$-cocycle condition
with the first term $x \cdot c(y,z,w)$ omitted.
This is expected from Fig.~\ref{d31}.
Constant functions, for example, satisfy this condition.
$\Box$

Next we look at a coboundary condition.
 A $3$-coboundary is written as
$$\xi_1(x \otimes y \otimes z) = \sum_u c_u(x,y,z) u
 = d^{2,1}(\phi) (x \otimes y \otimes z) =
  z^{-1} \phi (x  \otimes  y)  z
+ \phi(y^{-1}xy  \otimes  z) - \phi(x \otimes  yz) . $$ If we
write $\displaystyle{\phi(x,y)=\sum_u h_u (x,y) u}$, then
\begin{eqnarray*}
\lefteqn{(d^{2,1}(\phi) ) (x \otimes y \otimes z) }\\
&=&
  z^{-1} \left(\sum_u h_u (x,y) u \right) z
+\left( \sum_v  h_v (y^{-1}xy  ,z) v \right) - \left( \sum_w  h_w (x, yz) w \right)  \\
&=& \sum_g (\ h_{z g z^{-1} } (x, y) + h_{g} (y^{-1}xy  ,z) - h_g (x, yz) \ ) \ g .
\end{eqnarray*}
Hence
$$c_u (x,y,z) = h_{z u z^{-1} } (x, y) + h_{u} (y^{-1}xy  ,z) - h_u (x, yz) $$
and in particular for the coefficients $c_u(x,y,z)$ from
Proposition \ref{Gp3cocyProp},
$$c(x,y,z) =c_{z^{-1}y^{-1}xyz}  (x,y,z)
= h_{ y^{-1}xy} (x, y) + h_{z^{-1}y^{-1}xyz} (y^{-1}xy  ,z) -
h_{z^{-1}y^{-1}xyz} (x, yz) .$$ By setting $h_{ y^{-1}xy} (x,
y)=a(x,y)$, we obtain:
\begin{lemma}\label{G3coblem}
A $3$-cocycle $c(x,y,z)$ is
 a coboundary if for some
$a(x,y)$,
$$c(x,y,z) = a(x, y) + a (y^{-1}xy  ,z) - a(x, yz) .$$
\end{lemma}

\begin{remark}\label{Gpoidrem}
{\rm {}From Remark~\ref{G2rem}, Proposition~\ref{Gp3cocyProp}, and
Lemma~\ref{G3coblem}, we have the following situation. The
$2$-cocycle condition, the $3$-cocycle condition, and the
$3$-coboundary condition, respectively, gives rise to the
equations
\begin{eqnarray*}
& &   a(x, y) + a (y^{-1}xy  ,z) - a(x, yz) =0 ,\\
& & c(x, y, z)+ c(x, yz, w)-c(y^{-1} x y,  z, w) - c(x, y, zw) = 0, \\
& & c(x,y,z) = a(x, y) + a (y^{-1}xy  ,z) - a(x, yz) .
\end{eqnarray*}
This suggests a cohomology theory, which we investigate in Section~\ref{gpoidsec}.
} \end{remark}

\begin{proposition}
For the symmetric  group $G=S_3$ on three letters, we have
$H^1_{\ad}(kG;kG)=0$ and $H^2_{\ad}(kG;kG)\cong \displaystyle{\bigoplus_{3} (kG)}$
for $k=\C$
and $\F_3$.
\end{proposition}
{\it Proof.\/} By Lemma~\ref{kGC1lem}, we have
$H^1_{\ad}(kG;kG)=0$ and $B^2_{\ad}(kG,kG)=0$. Hence
$H^2(kG;kG)\cong Z^2_{\ad}(kG;kG)$, which is computed by solving
the system of equations stated in Lemma~\ref{G2lem} and
Remark~\ref{G2rem}. Computations by {\it Maple}
and {\it Mathematica}
shows that the
solution set is of dimension $3$ and generated by $(a((1\ 2\ 3),
(1\ 2))$, $a((2\ 3), (1\ 3\ 2) )$, and $a((1\ 3), (1\ 2))$ for the above
mentioned coefficient fields. $\Box$

\subsection{Function Algebras on Groups}

Let $G$ be a finite group and $k$ a field with char($k$) $\neq 2$.
The  set $k^G$ of functions from $G$ to $k$ with pointwise
addition and multiplication is a unital associative algebra.
It has a Hopf algebra structure using $k^{G \times G}\cong k^G \otimes k^G$
 with comultiplication
 defined through $\Delta: k^G \rightarrow k^{G \times G}$ by
 $\Delta(f)(u \otimes v)=f(uv)$  and the antipode by $S(f)(x)=f(x^{-1})$.

Now $k^G$ has basis (the characteristic function) $\delta_g : G
\rightarrow k$ defined by $\delta_g(x)=1$ if $x=g $ and zero
otherwise. Since   $S(\delta _g)=\delta_{g^{-1}}$ and
$\displaystyle{\Delta(\delta_h)=\sum_{uv=h} \delta_u \otimes
\delta_v}$,
  the adjoint map becomes
$$\ad(\delta_g \otimes \delta_h)
=\sum_{uv=h}
\delta_{u^{-1}} \delta_g  \delta_v =
\left\{ \begin{array}{ll} \delta_g \ & \mbox{\rm if $h=1$}, \\
                0\ & \mbox{\rm otherwise} .
                \end{array}
\right. $$

\begin{lemma}
$C^1_{\ad}(k^G; k^G)=0$.
\end{lemma}
{\it Proof.\/}
Recall that
$$C^1_{\rm ad}(H;H)=  \{ f \in \Hom_k (H,H) \ | \
f \mu = \mu (f \otimes 1) + \mu (1 \otimes f), \ \Delta f = (f
\otimes 1 ) \Delta + (1 \otimes f) \Delta \ \}.  $$ Let $G=\{ g_1,
\ldots, g_n \}$ be a given finite group and abbreviate
$\delta_{g_i}=\delta_i$ for $i=1, \ldots, n$. Describe  $f : k^G
\rightarrow k^G$ by $\displaystyle{f(\delta_i)=\sum_{j=1}^n
s_{i}^{ j} \delta_j}$. Then $f \mu = \mu (f \otimes 1) + \mu (1
\otimes f)$ is written for basis elements by ${\rm LHS} = f(\delta_i
\delta_j) $ and
\begin{eqnarray*}
{\rm RHS} &=& f(\delta_i) \delta_j + \delta_i f(\delta_j) \\
 &=& (\sum_{\ell =1}^n s_{i }^{ \ell} \delta_\ell ) \delta_j +
 \delta_i ( \sum_{h =1}^n s_{j}^{ h} \delta_h ) \\
 &=& s_{i}^{ j} \delta_j  + s_{j}^{ i } \delta_i.
 \end{eqnarray*}
 For $i=j$ we obtain ${\rm LHS} \displaystyle{= \sum_{w =1}^n s_{i }^{w} \delta_w}$
 and ${\rm RHS}=2 s_{i}^{ i} \delta_i$ so that
 $s_{i}^{ j}=0$ for all $i, j$ as desired.
 $\Box$

\begin{lemma}
$Z^2_{\ad}(k^G; k^G)=0$.
\end{lemma}
{\it Proof.\/} Recall that $d^{2,1}(\eta_1) = \ad(\eta_1\otimes
1)+ \eta_1 (\ad \otimes 1)- \eta_1(1 \otimes \mu) $ for $\eta_1
\in C^2_{\ad}(k^G, k^G)$. \vspace{2mm} Describe a general element
$\eta_1 \in C^2_{\ad}(k^G, k^G)$ by $\displaystyle{\eta_1(\delta_i
\otimes \delta_j )=\sum_{\ell} s_{i\,  j}^\ell\delta_\ell}$.
\vspace{2mm} Consider  $d^{2,1}(\eta_1)(\delta_a \otimes \delta_b
\otimes \delta_c )$. If $c\neq 1$, then the first term is zero by
the definition of $\ad$. If $c\neq 1 $ and $b=1$, then the third
term is also zero, and we obtain that the second term
$\eta_1(\delta_a \otimes \delta_c)$ is zero. Hence
$\eta_1(\delta_a \otimes \delta_c)=0$ unless $c=1$. Next, set
$b=c=1$ in the general form. Then all three terms equal
$\eta_1(\delta_a \otimes \delta_1)$ and we obtain $\eta_1(\delta_a
\otimes \delta_1)=0$, and the result follows. $\Box$

By combining the above lemmas, we obtain the following:

\begin{theorem}
For any finite group $G$ and a field $k$,
we have $H^n_{\ad}(k^G; k^G)=0$ for $n=1,2$.
\end{theorem}

Observe that $k(G)$ and $k^G$ are cohomologically distinct.

\subsection{Bosonization of the Superline}\label{bozonsubsec}

Let $H$ be generated by $1$, $g$, $x$ with relations
$x^2=0$, $g^2=1$, $xg=-gx$ and Hopf algebra structure
$\Delta(x)=x \otimes 1+g \otimes x$, $\Delta(g)=g \otimes g$, $\epsilon(x)=0$,
 $\epsilon(g)=1$, $S(x)=-gx$, $S(g)=g$
 (this Hopf algebra is called
 the bosonization of the superline~\cite{MajidPink}, page $39$, Example 2.1.7).

The operation $\ad$ is represented by the following table,
where, for example,  $\ad( g \otimes x)=2x$.

\bigskip

\begin{center}

\begin{tabular}{|r||r|r|r|r|}
\hline
       & $1$ & $g$ & $x$ & $gx$ \\   \hline \hline
$1$ & $1$ & $1$ & $0$ & $0$ \\  \hline
$g$ & $g$ & $g$ & $2x$& $2x$  \\  \hline
$x$ & $x$ & $-x$& $0$ & $0 $ \\  \hline
$gx$ & $gx$& $-gx$& $0$ & $0$ \\
\hline
\end{tabular}

\end{center}

\bigskip

\begin{remark}{\rm
The induced $R$-matrix $R_\ad$ has determinant $1$,
the characteristic polynomial
is
$(\lambda^2+1)^2 (\lambda+1)^4 (\lambda-1)^8$,
and the minimal polynomial
is
$(\lambda^2+1) (\lambda+1) (\lambda-1)^2$.
} \end{remark}

\begin{proposition}
The first cohomology of $H$ is given by $H^1_{\ad}(H,H)\cong k$.
\end{proposition}
{\it Proof.\/} Recall that $1$-cochains are given by
$$C^1_{\rm ad}(H;H)=  \{ f \in \Hom_k (H,H) \ | \
f \mu = \mu (f \otimes 1) + \mu (1 \otimes f), \
\Delta f = (f \otimes 1 ) \Delta + (1 \otimes f) \Delta \ \}.  $$
Let $f \in C^1_{\rm ad}(H;H)$.
Assume that $f(x)=a+bx+cg+dxg$ and $f(g)=\alpha +\beta x+\gamma g+ \delta xg$
 where $a,b,c,d, \alpha, \beta, \gamma, \delta  \in k$.
 Applying $f$ to both sides of the equation $g^2=1$,
 one obtains $\alpha=\gamma=0$.
Similarly evaluating both sides of the equation
$\Delta f = (f \otimes 1 ) \Delta + (1 \otimes f) \Delta $ at $g$ gives $\beta=\delta=0$,
one obtains  that $f(g)=0$.
In a similar way, applying $f$ to the equations $x^2=0$ and $xg=-gx$ gives
rise to, respectively,  $a=0$ and $c=0$.
Also evaluating $\Delta f = (f \otimes 1 ) \Delta + (1 \otimes f) \Delta $ at $x$
 gives rise to  $d=0$.
 We also have $f(x)=f(xg)g$ (since $g^2=1$),which implies
  $f(xg)=bxg$.  In conclusion $f$ satisfies $f(1)=0=f(g),f(x)=bx \;$,  and $f(xg)=b(xg)$.
Now consider $f$ in the kernel of $D_1$, that is $f$ satisfies
$$d^{1,1}(f)= \ad (1 \otimes f) - f\ad + \ad(f \otimes 1). $$
It is directly checked on all the generators $u \otimes v$ of $H \otimes H$ that $d^{1,1}(f)(u \otimes v)=0$.  This implies that $H^1(H,H)\cong k$.
$\Box$

\begin{proposition}\label{bozoprop}
For any field $k$ of characteristic not $2$,
 $H^2_{\ad}(H,H)\cong k^3$.
 \end{proposition}
 {\it Proof.\/}
 With $d^{1,1}=0$ from the preceding Proposition,
 we have $H^2_{\ad}(H,H)\cong Z^2_{\ad}(H,H)$.

For the convenience of the reader we compute,
 $\Delta(gx) = gx \otimes g + 1 \otimes gx$.
A number of key facts will be repeatedly recalled; these are inclosed in boxes.

The first $2$-differential is written as
$$\ad(\phi(a \otimes  b)\otimes c )+ \phi(\ad (a\otimes b )\otimes  c)
- \phi(a \otimes  bc)=0. $$

Take $b=c=1$, then since $\ad(a \otimes 1)=a$ for any $a \in H$,
all three
terms are the same and gives that \framebox{$\phi(a \otimes 1)=0$}
 for any $a$.

Take $a=g$ and  $b=c=x$, then the third term vanishes and we obtain
$\ad(\phi(g \otimes  x)\otimes x)+ \phi( 2x \otimes  x) =0$.
For any possible value of $\phi(g \otimes  x)$, the value of
the first term is written as $hx$ for some $h \in k$ from the table of $\ad$ above.
Since $\phi$ is bilinear,
constants can be renamed to obtain
\framebox{$\phi(x \otimes x)=hx$}.
A similar argument gives
\framebox{$\phi(x \otimes gx)=h'x$}
 from $a=g$, $b=x$ and  $c=gx$,
for another constant $h' \in k$.

The second differential is written as
$$\phi(a_{(1)} \otimes b_{(1)}) \otimes a_{(2)} b_{(2)} =
\phi(a \otimes b_{(2)} )_{(1)} \otimes b_{(1)} \phi(a \otimes b_{(2)} )_{(2)}. $$
Taking $a=b=x$, we obtain the LHS
$\phi (x \otimes x) \otimes 1 + (\phi (x \otimes g) + \phi (g \otimes x) )\otimes x.$
The RHS is $\phi(x \otimes x )_{(1)} \otimes g \phi(x \otimes x )_{(2)}$,
and using that $\phi(x \otimes x )=hx$, we obtain
$h(x \otimes g + g \otimes gx)$ for the RHS.
Since there is no $\otimes 1$ term in the RHS, we obtain $\phi (x \otimes x)=0$,
and in particular, $h=0$, which makes ${\rm RHS}=0$, and we also obtain
\framebox{$\phi(x \otimes g)= - \phi(g \otimes x)$}.

Let $a=b=g$ in the second differential.
Then ${\rm LHS} =\phi(g \otimes g) \otimes 1$ and
${\rm RHS} =\phi(g \otimes g)_{(1)} \otimes g \phi(g \otimes g)_{(2)}$.
This implies that $\phi(g \otimes g)$ is written as $h_g g$ for some $h_g \in k$,
and ${\rm RHS} =
h_g
(g \otimes g^2 )= h_g  (g \otimes 1)={\rm LHS.}$
With $a=b=c=g$ in the first differential, we obtain
$\ad(h_g g \otimes g)+ h_g g - 0 = 0$, hence, in fact, $h_g=0$ if $2$ is invertible,
giving rise to
\framebox{$\phi(g \otimes g)=0$}.

Let $a=g$ and $b=x$  in the second differential.
Then the
 ${\rm LHS} =\phi(g \otimes x)\otimes g$,
and
the
${\rm RHS} =\phi(g \otimes x)_{(1)} \otimes g \phi(g \otimes x)_{(2)} $.
{}For the
RHS to have terms
ending in
$\otimes g$
only,
$\phi(g \otimes x)$
can have neither
$g$ nor $gx$ terms since they would
result in a $(\ \otimes 1)$
term, so let
$\phi(g \otimes x)= h_{g,x} 1 + h_{g,x}' x$.
Then one computes
${\rm RHS} =( h_{g,x} 1 +  h_{g,x}' x) \otimes g + h_{g,x}' (g \otimes gx)$.
Equating this with LHS, we obtain $h_{g,x}'=0 $.
Thus we obtained
\framebox{$\phi(g \otimes x)= h_{g,x} 1 =- \phi(x \otimes g)$}. 
In the first differential, take $a=b=g$ and $c=x$ to obtain
\framebox{$\phi(g \otimes gx)= \phi(g \otimes x)= h_{g,x} 1$}.

Let $a=1$ and $b=x$  in the second differential.
Then the
 ${\rm LHS} =\phi(1 \otimes x)\otimes 1 + \phi(1 \otimes g)\otimes x$,
and
the
${\rm RHS} =\phi(1 \otimes x)_{(1)} \otimes g \phi(1 \otimes x)_{(2)} $.
{}For the
RHS to have terms
ending in
$\otimes 1$ or $\otimes x$
only,
$\phi(1 \otimes x)$
can have neither
$1$ nor $x$ terms since they would
result in a $(\ \otimes g)$
term, so let
\framebox{$\phi(1 \otimes x)= h_{1, x} g + h_{1, g} gx $}.
Then one computes
${\rm RHS} = h_{1,x} g\otimes 1  +  h_{1, g} (gx \otimes 1 + 1 \otimes x)$.
Comparing with the LHS, we obtain
\framebox{$\phi(1 \otimes g)=h_{1, g} 1$}.
With $a=1$, $b=x$ and $c=g$ in the first differential, we also obtain
\framebox{$\phi(1 \otimes gx)=-h_{1, x} g + h_{1, g} gx $}.

Recall that $\phi(x \otimes gx)=h'x$.
 For $a=x$ and $b=gx$ in the second differential gives
 \begin{eqnarray*}
 {\rm LHS}&=&\phi( x \otimes gx) \otimes g - \phi(g \otimes gx) \otimes gx
 =h' (x  \otimes g) - h_{g,x} (1 \otimes gx) \\
 {\rm  RHS} &=& - h_{g, x} (1 \otimes gx)
  +h' (x \otimes 1 + g \otimes x)
 \end{eqnarray*}
 which implies
 \framebox{$ \phi( x \otimes gx) = 0$}.

In the second differential,
take $a=gx$ and $b=x$.
Then we obtain
 \begin{eqnarray*}
 {\rm LHS}&=&\phi( gx \otimes x) \otimes g +
 (\phi(gx \otimes g) + \phi(1 \otimes x) ) \otimes gx \\
 &= &   \phi( gx \otimes x) \otimes g +
 (\phi(gx \otimes g) + h_{1, x} g + h_{1, g} gx   ) \otimes gx\\
 {\rm  RHS} &=&
   \phi(gx \otimes x)_{(1)} \otimes g \phi(gx \otimes x)_{(2)} .
 \end{eqnarray*}
The LHS has only $\otimes g $ and $\otimes gx$ terms, so
that $\phi(gx \otimes x)$ does not have $g$ or $gx$ terms,
and we can write $\phi(gx \otimes x)= h_{gx,x} 1 +  h_{gx,x}' x$ and compute
${\rm RHS} =h_{gx,x} (1\otimes g) +  h_{gx,x}' (x \otimes g + g \otimes gx )$.
Comparing with the LHS we obtain
 $ h_{gx,x}' g=\phi(gx \otimes g) + h_{1, x} g + h_{1, g} gx  $,
 so that \framebox{$\phi(gx \otimes g)=  ( h_{gx,x}' - h_{1, x} ) g - h_{1, g} gx  $}.

 By the first differential with $(a,b,c)=(gx, g, x)$,
 we obtain
 $$\phi(gx \otimes gx)= 2 ( h_{gx,x}' - h_{1, x} ) x - (h_{gx,x} 1 +  h_{gx,x}' x)
 = - h_{gx,x} 1 + ( h_{gx,x}' - 2 h_{1, x} ) x. $$

\begin{sloppypar}
 By the first differential with $(a,b,c)=(gx, gx, g)$,
 we obtain
$$ ( - h_{gx,x} 1  -
( h_{gx,x}' - 2 h_{1, x} ) x ) + 0 + (h_{gx,x} 1 +  h_{gx,x}' x)=0 $$
which implies $h_{1,x}=0$.
In particular,
we obtain \framebox{$\phi(gx \otimes g)=h_{gx,x}' g - h_{1, g} gx$}
and
\framebox{$\phi(gx \otimes gx)=  - h_{gx,x} 1 +  h_{gx,x}'  x$}.
By the second differential with $a=b=gx$,
 we obtain
\begin{eqnarray*}
 {\rm LHS}&=&
 \phi(gx \otimes gx) \otimes 1 + \phi(1 \otimes gx) \otimes (gx)g \\
 &=&( -h_{gx,x} 1 + h_{gx,x}' x ) \otimes 1 - h_{1,g} ( gx  \otimes x) \\
  {\rm  RHS} &=&
  \phi(gx \otimes g)_{(1)} \otimes (gx)   \phi(gx \otimes g)_{(2)}
  +   \phi(gx \otimes gx)_{(1)} \otimes    \phi(gx \otimes gx)_{(2)}\\
  &=&
(h_{gx,x}' (g \otimes gxg) + h_{1,g} (gx \otimes x) )
+ (- h_{gx,x} (1 \otimes 1) + h_{gx,x}' (x \otimes 1 + g \otimes x))
 \end{eqnarray*}
and comparing the terms we obtain
$2 h_{1,g}=0$.
In summary, resetting free variables by
$h_{g,x}=\alpha$, $h_{gx,x}=\beta$ and $h_{gx,x}'=\gamma$,
we obtained a general solution represented by the following table.
\end{sloppypar}

\bigskip

\begin{center}

\begin{tabular}{|r||r|r|r|r|}
\hline
       & $1$ & $g$ & $x$ & $gx$ \\   \hline \hline
$1$ & $0$ & $0$ & $0$ & $0$ \\  \hline
$g$ & $0$ & $0$ & $\alpha 1 $& $ \alpha 1$  \\  \hline
$x$ & $0$ & $- \alpha  1$& $0$ & $0 $ \\  \hline
$gx$ & $0$& $\gamma g $& $\beta  1 + \gamma x $ & $- \beta  1 + \gamma x$ \\
\hline
\end{tabular}

\end{center}

\bigskip

\noindent
It is checked, either by hand, or computer guided calculations, that
these 
are indeed solutions.
$\Box$

\section{Adjoint, Groupoid, and Quandle Cohomology Theories}\label{gpoidsec}

{}From Remark~\ref{Gpoidrem},
the adjoint cohomology leads us to
cohomology, especially for
conjugate groupoids of groups as defined below.
Through the relation between Reidemeister moves for knots
and the adjoint, groupoid cohomology, we obtain a new construction of
quandle cocycles. In this section we investigate these relations.
First we formulate a general definition.
Many formulations of groupoid cohomology can be found  in literature,
and relations of the following formulation to previously  known theories are not clear.
See \cite{Tu}, for example.

Let $\GG$ be a groupoid with objects $\textrm{Ob}({\GG})$ and
morphisms $G(x,y)$ for $x, y \in \textrm{Ob}(\GG)$. Let $f_i \in
G(x_i, x_{i+1})$, $0 \leq i < n$,
 for non-negative integers $i$ and $n$.
Let $C_n(\GG)$ be the free abelian group generated by
$$ \{ (x_0, f_0, \ldots, f_n) \ | \ x_0 \in \textrm{Ob}(\GG), f_i \in G(x_i, x_{i+1}), 0 \leq i < n \}. $$
The boundary map
$\partial : C_{n+1}(\GG) \rightarrow C_{n}(\GG) $
is defined by
by linearly extending
\begin{eqnarray*}
\lefteqn{\partial (x_0, f_0, \ldots, f_n)= (x_1, f_1, \ldots, f_n) }\\
&+ &\sum_{i=0}^{n-1} (-1)^{i+1} (x_0, f_0, \ldots, f_{i-1},  f_i f_{i+1}, f_{i+2}, \ldots, f_n)\\
&+ &(-1)^{n+1} (x_0, f_0, \ldots, f_{n-1}) .
\end{eqnarray*}
Then it is easily seen that this differential defines a chain complex.

The corresponding groupoid 1- and 2-cocycle conditions are written as:
\begin{eqnarray*}
& & a(x_1, f_1) - a(x_0, f_0 f_1) + a(x_0, f_0) = 0 \\
& & c(x_1, 
f_1, f_2) - c(x_0, f_0 f_1, f_2) + c(x_0, f_0, f_1 f_2) - c(x_0, f_0, f_1) = 0
\end{eqnarray*}  The general cohomological theory of homomorphisms and extensions applies, such as:
\begin{remark}{\rm Let $\GG$ be a groupoid and $A$ be an abelian group regarded as a one-object groupoid.  Then $\alpha: \hom(x_0, x_1) \rightarrow A$ gives a groupoid homomorphism from $\GG$ to $A$, which sends $\rm{Ob}(\GG)$ to the single object of $A$, if and only if $a: C_1(\GG) \rightarrow A$, defined by $a(x_0, f_0) = \alpha(f_0)$, is a groupoid 1-cocycle.

Next we consider extensions of groupoids.  Define $\circ: (\hom(x_0, x_1) \times A) \times (\hom(x_1, x_2) \times A) \rightarrow \hom(x_0, x_2) \times A$ by $$(f_0, a) \circ (f_1, b) = (f_0 f_1, a+b+c(x_0, f_0, f_1))$$ where $c(x_0, f_0, f_1) \in \hom(C_2 (\GG), A)$.  If $\GG \times A$ is a groupoid,  
the function $c$ with the value 
 $c(x_0, f_0, f_1)$ is a groupoid 2-cocycle.
}
\end{remark}

\begin{figure}[htb]
\begin{center}
\mbox{
\epsfxsize=3in
\epsfbox{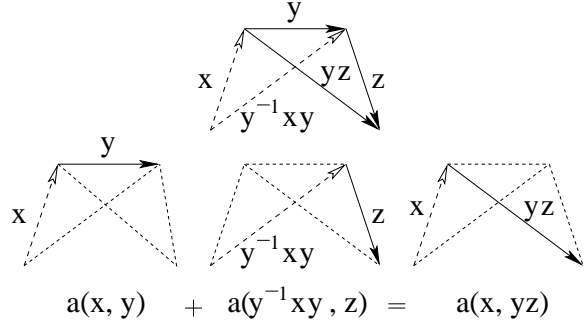}
}
\end{center}
\caption{Diagrams for a groupoid $1$-cocycle}
\label{gpoid1c}
\end{figure}

\begin{figure}[htb]
\begin{center}
\mbox{
\epsfxsize=3in
\epsfbox{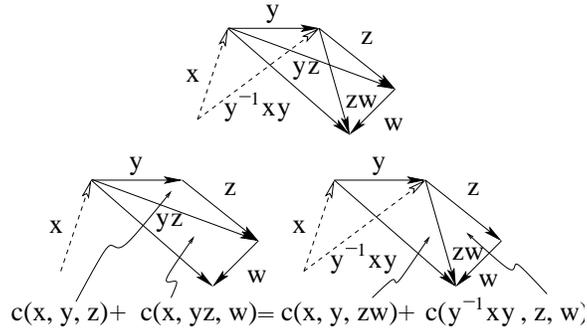}
}
\end{center}
\caption{Diagrams for a groupoid $2$-cocycle}
\label{gpoid2c}
\end{figure}

\begin{example}{\rm
Let $G$ be a group. Define the
{\it conjugate groupoid} of $G$, denoted $\widehat{G}$, by:
\begin{eqnarray*}
\textrm{Ob}(\widehat{G}) &=& G \\
\textrm{Mor}(\widehat{G}) &=& G \times G
\end{eqnarray*}
where the source of the morphism $(x,y) \in \hom(x, y^{-1}xy)$ is $x$ and its target is $y^{-1}xy$, for $x,y \in G$.  Composition is defined by $(x,y) \circ (y^{-1}xy, z) = (x, yz)$.  For this example, the groupoid 1- and 2-cocycle conditions are:
\begin{eqnarray*}
& &   a(x, y) + a (y^{-1}xy  ,z) - a(x, yz) =0 ,\\
& & c(x, y, z)+ c(x, yz, w)-c(y^{-1} x y,  z, w) - c(x, y, zw) = 0.
\end{eqnarray*}
Diagrammatic representations of these equations are
depicted in Figs.~\ref{gpoid1c}, \ref{gpoid2c}.
Furthermore, $c$ is a coboundary if $$c(x,y,z) = a(x, y) + a (y^{-1}xy  ,z) - a(x, yz) .$$
Compare with Remark~\ref{Gpoidrem}.

For $G=\Sym_3$, the symmetric group on $3$ letters,
with coefficient group $\C$, $\Z_2$, $\Z_3$, $\Z_5$ and $\Z_7$,
respectively, the dimensions of the conjugation groupoid $2$-cocycles are
$3$, $5$, $4$, $3$ and $3$.

}\end{example}

For the rest of the section, we present new constructions of
quandle cocycles from groupoid cocycles of conjugate groupoids
of groups.
Let $G$ be a finite group, and $a: G^2 \rightarrow k$ be  adjoint $2$-cocycle
coefficients that were defined
in Remark~\ref{G2rem}.  These satisfy
$$a (x, y)+a(x \lt y, z)-a(x, yz)=0. $$

\begin{proposition}
Let $\psi(x, y)=a(x,y)$.
Then $\psi$ satisfies the rack $2$-cocycle condition
$$\psi(x, y)+\psi(x \lt y, z)=\psi(x, z) + \psi(x \lt z, y \lt z).$$
\end{proposition}
{\it Proof.\/}
By definition $\psi(x, y)+\psi(x \lt y, z) = a(x,yz)$, and $\psi(x, z) + \psi(x \lt z, y \lt z) =
a(x, z) + a(z^{-1}x z, z^{-1}y z)=a(x , z(y\lt z))$.
$\Box$

Let $G$ be a finite group, and $c: G^3 \rightarrow k$ be a coefficient of the adjoint $3$-cocycle defined
in Proposition~\ref{Gp3cocyProp}.   This satisfies
$$c(x, y, z)+ c(x, yz, w)=c(x \lt y, z, w) + c(x, y, zw).$$
 \begin{proposition}
 Let $G$ be a group that is considered as a quandle under conjugation.
 Then  $\theta: G^3 \rightarrow k$ defined by $\theta(x,y,z)=c(x,y,z)-c(x, z, z^{-1}yz)$ is a rack $3$-cocycle. \end{proposition}
 {\it Proof.\/}
 We must show that $\theta$ satsifies
\begin{eqnarray*}
\lefteqn{ \theta(x,y,z) + \theta(x\lt z, y \lt z, w) + \theta(x, z, w) }\\
&=&
 \theta(x\lt y,z, w) + \theta(x, y, w) + \theta(x\lt w, y \lt w, z\lt w) .
 \end{eqnarray*}
 We compute
\begin{eqnarray*}
{\rm LHS } - {\rm  RHS} &=&
[ \ c(x,y,z) - c(x, z, z^{-1}yz) \ ] \\
& & + [\ c(z^{-1}xz, z^{-1}yz, w) - c(z^{-1}xz, w, w^{-1} z^{-1}yz w )\ ]\\
& & + [ \ c(x, z, w) - c(x, w, w^{-1}zw) \ ] \\
& & -
[ \ c(y^{-1}xy, z, w) - c(y^{-1}xy, w, w^{-1}zw)\ ] \\
& & - [\ c(x, y, w) - c(x, w, w^{-1}yw) \ ] \\
& & - [\ c(w^{-1} x w, w^{-1}y  w, w^{-1}z w) - c(w^{-1} x w, w^{-1}z  w, w^{-1} z^{-1}y zw) \ ] \\
& = & [ \ c(x,y,z) -  c(y^{-1}xy, z, w) \ ] \\
& & -  [\  c(x, z, z^{-1}yz) - c(z^{-1}xz, z^{-1}yz, w) \ ] \\
& & + [ \ c(x, z, w) -  c(z^{-1}xz, w, w^{-1} z^{-1}yz w )\ ]\\
& & - [ \  c(x, w, w^{-1}zw) - c(w^{-1} x w, w^{-1}z  w, w^{-1}z^{-1}y zw) \ ] \\
& & - [\ c(x, y, w) -  c(y^{-1}xy, w, w^{-1}zw)\ ] \\
& & + [ \  c(x, w, w^{-1}yw) - c(w^{-1} x w, w^{-1}y  w, w^{-1}z w)\ ] \\
& = & [\ - c(x, yz, w) + c(x,y,zw) \ ] \\
& & - [ \ - c(x, z z^{-1} yz, w) + c(x,z,z^{-1}yzw) \ ] \\
& & + [\ - c(x, zw, w^{-1} z^{-1} yzw) + c(x,z,w w^{-1} z^{-1} yzw) \ ] \\
& & - [ \ - c(x, w  w^{-1}zw ,  w^{-1}z^{-1}yzw) + c(x,w,w^{-1}z w w^{-1}y w)\ ] \\
& & -[ \ -c(x, yw, w^{-1} z w ) + c(x,y,w w^{-1} z w ) \ ] \\
& & + [ \ - c(x, w w^{-1}yw, w^{-1} z w)
+ c( x,w, w^{-1}y  w w^{-1}z w)\ ] \\
& = & \ 0
\end{eqnarray*}
as desired. $\Box$

\section{Deformations of $R$-matrices by adjoint $2$-cocycles}
\label{Rmat}

In this section we give, in an explicit form, deformations of
R-matrices by $2$-cocycles of the adjoint cohomology theory we
developed in this paper. Let $H$ be a Hopf algebra and $\ad$ its
adjoint map. In Section~\ref{deformsec} a deformation of $(H, \ad)
$ was defined to be a pair $(H_t, \ad_t)$ where $H_t$ is a
$k[[t]]$-Hopf algebra given by $H_t=H \otimes k[[ t ]]$ with all
Hopf algebra structures inherited by extending those on $H_t$. Let
$A=(H \otimes  k[[ t ]] ) / (t^2) )$
and the Hopf algebra structure
maps  $\mu, \Delta, \epsilon, \eta, S$ be inherited on $A$. As a
vector space $A$ can be regarded as $H \oplus tH$

Recall that a solution to the YBE, $R$-matrix $R_{\ad}$ is induced from
the adjoint map.
Then from the constructions of the adjoint cohomology from
the point of view of the deformation theory, we obtain the following deformation
of this $R$-matrix induced from the adjoint map.

\begin{theorem}
Let $\phi \in Z^2_\ad (H;H)$ be an adjoint $2$-cocycle.
Then the  map $R: A \otimes A \rightarrow A \otimes A$ defined by
$R=R_{\ad + t \phi}$ satisfies the YBE.
\end{theorem}
{\it Proof.\/}
The equalities of Lemma~\ref{deformlem}
hold
 in
 the
  quotient
$A=(H \otimes  k[[ t ]]) / (t^2) $, 
where $n=1$ and the modulus  
$t^2$ is considered.
These cocycle conditions, on the other hand,
were formulated from the motivation  from Lemma~\ref{adYBElem}
for the induced $R$-matrix $R_\ad$ to satisfy the YBE.
Hence these two lemmas imply the theorem.
$\Box$

\begin{example}{\rm
In Subsection~\ref{bozonsubsec},
the adjoint map $\ad$ was computed
for the bosonization $H$ of the superline,
with basis $\{1, g, x, gx\}$,
as well as a general $2$-cocycle $\phi$
with three free variables $\alpha, \beta, \gamma$
written by $\phi(g \otimes x)=\phi(g \otimes gx)=\alpha 1$,
$\phi(x \otimes g)= - \alpha 1$,
$\phi(gx \otimes g)=\gamma g$,
$\phi(gx \otimes x)=\beta 1 + \gamma x$,
$\phi(gx \otimes gx)= - \beta 1 + \gamma x$,
and zero otherwise.
Thus we obtain the deformed solution to the YBE $R=R_{\ad + t \phi}$
on $A$ with three variables $t \alpha, t \beta, t \gamma$ of degree one.
} \end{example}

\section{Concluding Remarks}

In \cite{CCES1} we concluded with  {\it A Compendium of Questions}
regarding our discoveries. Here we attempt to address some of
these questions by providing relationships between this paper and
\cite{CCES1}, and offer further questions for our future
consideration.

It was pointed out in  \cite{CCES1}
 that there was a clear distinction between
the Hopf 
algebra case and the cocommutative coalgebra case as to why
self-adjoint maps satisfy the YBE. In  \cite{CCES1} a cohomology
theory was constructed for the coalgebra case. In this paper, many
of the same ideas and techniques, in particular deformations and
diagrams, were used to construct a cohomology theory in the Hopf
algebra case, with applications to the YBE and quandle cohomology.

The aspects that unify these  two theories are deformations and
a systematic process we call ``diagrammatic infiltration."
So far,
these techniques have only been successful in defining
coboundaries up through dimension $3$. This is a deficit of the
diagrammatic approach, but
diagrams give direct applications to other algebraic problems
such as the YBE and quandle cohomology, and suggest further applications to
knot theory. By taking the trace as in Turaev's~\cite{Tur},
for example, a new deformed version of a given invariant is expected to
be obtained.

Many questions remain:
Can $3$-cocycles be used for solving the tetrahedral equation?
Can they be used for knotted surface invariants?
 Can the
coboundary maps be expressed skein theoretically?
How are the deformations of $R$-matrices related to deformations of underlying
Hopf algebras?
When a Hopf algebra contains a coalgebra, such
as the universal enveloping algebra and its Lie algebra together with the
ground field of degree-zero part, what is the relation between the two theories
developed in this paper and in  \cite{CCES1}?
How these theories, other than the same diagrammatic techniques,
can be uniformly formulated, and to higher dimensions?


\begin{thebibliography}{99}
\setlength{\itemsep}{-2pt}


\bibitem{AG}
{Andruskiewitsch, N.; Gra\~{n}a, M.,}
{\it From racks to pointed Hopf algebras},
Adv. in Math. {\bf 178} (2003), 177--243.

\bibitem{BC} Baez, J.C.; Crans, A.S.,
{\it Higher-Dimensional Algebra VI: Lie 2-Algebras}, Theory and
Applications of Categories {\bf 12} (2004), 492--538.

\bibitem{BL98}
{Baez, J.C.;   Langford, L.,}
{\it $2$-tangles},
{Lett. Math. Phys.}
{\bf 43} (1998), 187--197.

\bibitem{Br88} Brieskorn, E.,
{\it Automorphic sets and singularities,}
Contemporary math., {\bf 78} (1988), 45--115.

\bibitem{CES}
{J. S. Carter; M. Elhamdadi;  M. Saito}
{\it Twisted Quandle homology theory and cocycle knot invariants}
{Algebraic and Geometric Topology}
{\bf 2} (2002), 95--135.

\bibitem{CJKLS}
 Carter, J.S.; Jelsovsky, D.; Kamada, S.; Langford, L.; Saito, M.,
{\it Quandle cohomology and state-sum invariants
of knotted curves and surfaces,}
Trans. Amer. Math. Soc.  {\bf 355}  (2003),  3947--3989.

\bibitem{CCES1}
Carter, J.S.; Crans, A.; Elhamdadi, M.; Saito, S., {\it Cohomology
of Categorical Self-Distributivity}, Preprint,  available at
arXiv:math.GT/0607417.

\bibitem{Alissa} Crans, A.S.,
{\it Lie $2$-algebras}, Ph.D. Dissertation, 2004, UC Riverside,
available at arXiv:math.QA/0409602.

\bibitem{FR}   Fenn, R.; Rourke,  C.,
\textit{Racks and links in codimension two,}
Journal of Knot Theory and Its Ramifications {\bf  1}  (1992), 343--406.

\bibitem{Gerst} Gerstenharber, M; Schack, S.D.,
{\it Bialgebra cohomology, deformations, and quantum groups,}
Proc. Nat. Acad. Sci. U.S.A., {\bf 87} (1990), 478--481.

\bibitem{Henn}
Hennings, M.A.,
{\it On solutions to the braid equation identified by Woronowicz,}
Lett. Math. Phys. {\bf 27} (1993), 13--17.

\bibitem{Joyce} Joyce, D.,
{\it A classifying invariant of knots, the knot quandle,}
J. Pure Appl. Alg. {\bf 23}
(1982) 37--65.


\bibitem{Kuperberg} Kuperberg, G., {\it Involutory Hopf algebras and $3$-manifold invariants},  Internat. J. Math.  {\bf 2 } (1991),   41--66.

\bibitem{MajidGreen} Majid, S. ``A quantum groups primer.'' London Mathematical Society Lecture Note Series, 292. Cambridge University Press, Cambridge, 2002.


\bibitem{MajidPink} Majid, S. ``Foundations of quantum group theory.'' Cambridge University Press, Cambridge, 1995.


\bibitem{MrSt}
Markl, M.; Stasheff, J.D.,
{\it Deformation theory via deviations,}  J. Algebra  {\bf 170}  (1994),   122--155.


\bibitem{MrVo}  Markl, M.; Voronov, A.; {\it PROPped up graph cohomology},
to appear in Maninfest, preprint at http://arxiv.org/pdf/math/0307081.



\bibitem{Matveev} Matveev, S.,
{\it Distributive groupoids in knot theory,} (Russian) Mat. Sb. (N.S.)
{\bf 119(161)} (1982), 78--88 (160).


\bibitem{ReshTur} Reshetikhin, N.; Turaev, V. G., {\it Invariants of $3$-manifolds via link polynomials and quantum groups},  Invent. Math.  {\bf 103}  (1991),   547--597.


\bibitem{Tu}
Tu, J.-L.,
{\it Groupoid cohomology and extensions},
Trans. Amer. Math. Soc. {\bf 358} (2006), 4721--4747.

\bibitem{Tur} Turaev, V. G. {\it The Yang-Baxter equation and invariants of links},  Invent. Math.  {\bf 92}  (1988),  no. 3, 527--553.


\bibitem{Woro}
 Woronowicz, S.L.,
{\it   Solutions of the braid equation related to a Hopf algebra,}
  Lett. Math. Phys.  {\bf 23}  (1991),   143--145.


\end{thebibliography}
\end{document}